\documentclass{article}
\usepackage{latexsym}
\usepackage{amsfonts}
\usepackage{graphicx}

\newtheorem{theorem}{Theorem}[section]
\newtheorem{lemma}{Lemma}[section]

\newtheorem{definition}{Definition}[section]
\newtheorem{prop}{Proposition}[section]
\newcommand{\qed}{\hfill$\Box$\par\medskip}

\newenvironment{Proof}{\noindent{\sc Proof.}}{\qed}

\def\bhag#1{\noindent
\setcounter{equation}{0}
\section{#1}
}

\newcounter{udaharana}
\newenvironment{example}{\noindent\textsc{Example}  \refstepcounter{udaharana} \arabic{udaharana}.}{\qed}

 \renewcommand{\theequation}{\arabic{section}.\arabic{equation}}
\def\bfgk#1{{{#1}\kern-5.5pt{#1}}}

\def\RR{{\mathbb R}}
\def\CC{{\mathbb C}}
\def\XX{{\mathbb X}}
\def\ZZ{{\mathbb Z}}

\def\PPI{{{\rm I}\kern-1pt\Pi}}

\def\a{\alpha}
\def\b #1;{{\bf #1}}
\def\x{{\bf x}}
\def\k{{\bf k}}
\def\m{{\bf m}}
\def\y{{\bf y}}

\def\u{{\bf u}}
\def\e{\epsilon}

\def\C{{\mathcal C}}

\def\X{{\mathcal X}}

\def\Y{{\mathcal Y}}

\def\ip#1#2{{\langle {#1}, {#2}\rangle}}

\def\esssup{\mathop{\hbox{{\rm ess sup}}}}

\def\be{\begin{equation}}
\def\ee{\end{equation}}
\def\bea{\begin{eqnarray}}
\def\eea{\end{eqnarray}}
\def\eref#1{(\ref{#1})}
\def\disp{\displaystyle}

\def\span{\mbox{{\textsf{span }}}}
\def\donchitre#1#2{\vskip 6.5cm\noindent
\parbox[t]{1in}{\special{eps:#1.eps x=6.5cm y=5.5cm}}
\hbox to 7cm{}\parbox[t]{0.0cm}{\special{eps:#2.eps x=6.5cm y=5.5cm}}}

\def\tn{|\!|\!|}
\def\dist{\mbox{\textsf{ dist }}}
\def\supp{\mbox{\textsf{supp}}}

\begin{document}

\title{Marcinkiewicz--Zygmund measures on manifolds}         
\author{F.~Filbir\thanks{
Institute of Biomathematics
and Biometry,
Helmholtz Center Munich,
85764 Neuherberg, Germany,
\textsf{email:} filbir@helmholtz-muenchen.de. The research of this author was partially funded by Deutsche Forschungsgemeinschaft grant FI 883/3-1 and PO711/9-1.},
H.~N.~Mhaskar\thanks{Department of Mathematics, California State University,
Los Angeles, California, 90032, USA,
\textsf{email:} hmhaska@gmail.com. The
    research of this author was supported, in part, by grant  
    DMS-0908037 from the National Science Foundation and grant
    W911NF-09-1-0465 from the U.S. Army Research Office.}}        
\date{}          
\maketitle

\begin{abstract}
Let $\XX$ be a compact, connected, Riemannian manifold (without boundary), $\rho$ be the geodesic distance on $\XX$, $\mu$ be a probability measure on $\XX$, and $\{\phi_k\}$ be an orthonormal system of continuous functions, $\phi_0(x)=1$ for all $x\in\XX$, $\{\ell_k\}_{k=0}^\infty$ be an nondecreasing sequence of real numbers with $\ell_0=1$, $\ell_k\uparrow\infty$ as $k\to\infty$, $\Pi_L:=\span\{\phi_j : \ell_j\le L\}$,  $L\ge 0$. We describe conditions to ensure an equivalence between the $L^p$ norms of elements of $\Pi_L$ with their suitably discretized versions. We also give intrinsic criteria to determine if any system of weights and nodes allows such inequalities. The results are stated in a very general form, applicable for example, when the discretization of the integrals is based on weighted averages of the elements of $\Pi_L$ on geodesic balls rather than point evaluations.
\end{abstract}

\bhag{Introduction}
To avoid complicating our notations unnecessarily, the notations used in the introduction and the next section might have a different meaning from the rest of the paper. 

The classical Marcinkiewicz--Zygmund (MZ) inequality states the following \cite[Chapter~X, Theorem~(7.5)]{zygmund}:
Let $n\ge 1$ be an integer, and $S$ be a trigonometric polynomial of order at most $n$ (i.e., an expression of the form $\sum_{|k|\le n} c_ke^{ikx}$). If $1<p<\infty$, then
\be\label{mzorig}
\int_0^{2\pi}|S(x)|^pdx\le \frac{A_p}{2n+1}\sum_{k=0}^{2n}\left|S\left(\frac{2k\pi}{2n+1}\right)\right|^p\le AA_p\int_0^{2\pi}|S(x)|^pdx,
\ee
where $A$ is an absolute, positive constant, and $A_p$ is a positive constant depending only on $p$. We observe that the number of points in the summation is the same as the dimension of the space of all trigonometric polynomials of order at most $n$. The upper inequality in \eref{mzorig} is valid for  $p=1,\infty$ as well. The lower inequality holds  for $p=1,\infty$ if one allows more points in the summation than the dimension $2n+1$ (\cite[Chapter~X, Theorem~(7.28)]{zygmund}). These inequalities are also known as large sieve inequalities or network inequalities. Inequalities of this form have many applications in approximation theory, number theory, signal processing, etc. Therefore, several analogues of these inequalities have been studied in the literature, for example, in the setting where $S$ is a univariate algebraic polynomial, the integrals in \eref{mzorig} are replaced by weighted or Lebesgue--Stieltjes integrals on real intervals, and the weights  and sampling nodes in the sum in \eref{mzorig} are chosen judiciously. A survey of many of the classical results in this direction and their applications can be found in the paper \cite{lubinskymz} by Lubinsky. 

To describe the MZ inequality in a very general setting, let $\XX$ be a compact set, $\mu$ be a finite Borel measure on $\XX$, $\{\phi_k\}$ be an orthonormal system of continuous functions,  $\{\ell_k\}_{k=0}^\infty$ be an nondecreasing sequence of real numbers with $\ell_0=1$, $\ell_k\uparrow \infty$ as $k\to\infty$, $\Pi_L:=\span\{\phi_j : \ell_j\le L\}$,  $L\ge 0$. For each $L\ge 1$, let $\C_L$ be a finite subset of $\XX$, $W_L=\{w_y\}_{y\in\C_L}\subset \RR$. Let $1\le p\le \infty$. For the purpose of this introduction, we will say that $\{(\C_L,W_L)\}$ is a MZ system of order $p$ if for every $L\ge 1$ and $P\in\Pi_L$,
\be\label{mzgen1}
c_1\int_\XX |P(x)|^p d\mu(x) \le \sum_{y\in\C_L} |w_y||P(y)|^p \le c_2\int_\XX |P(x)|^p d\mu(x),
\ee
where $c_1,c_2$ are positive constants depending only on $\XX$, $\{\phi_k\}$, $\{\ell_k\}$, $\mu$, and $p$, but independent of $L$ and $P$, and the customary interpretation is assumed when $p=\infty$. We list four of the applications of such inequalities, which have inspired our own interest in these. 

First, suppose one wants to approximate a continuous function $f:\XX\to\RR$. A standard way to do this is  by means of an operator
$$
T_L(f,x)=\int_\XX f(y)\Phi_L(x,y)d\mu(y),
$$
where $\Phi_L$ is a suitable kernel with the property that it is in $\Pi_L$ as a function of $x$ and as a function of $y$. A typical example is the Fourier projection, where $\Phi_L(x,y)=\sum_{\ell_j\le L}\phi_j(x)\phi_j(y)$ is the reproducing kernel for $\Pi_L$. If the approximation is required using the values of $f$ at points in $\C_L$, then it is natural to consider the discretization
$$
T_L^D(f,x)=\sum_{y\in\C_L} w_y f(y)\Phi_L(x,y).
$$
In the case when $T_L$ is the Fourier projection, such a discretization has been called hyperinterpolation. It is easy to verify that the operator norms of $T_L$, $T_L^D$ are given by
$$
\sup_{x\in\XX}\int_\XX |\Phi_L(x,y)|d\mu(y), \qquad \sup_{x\in\XX}\sum_{y\in\C_L} |w_y| |\Phi_L(x,y)|
$$
respectively. Since $\Phi_L(x,\circ)\in \Pi_L$ for each $x\in\XX$, \eref{mzgen1} with $p=1$ implies that the operator norms of $T_L$, $T_L^D$ have the same order of magnitude as functions of $L$, as $L\to\infty$. In this context, it is not necessary that the weights $w_y$ should be all nonnegative, a fact which may be useful in numerical computations. In particular, in the case when $\XX$ is a Euclidean sphere, $\Phi_L$ is the reproducing kernel for the space $\Pi_L$ of spherical polynomials of degree at most $L$, then this leads to a simple proof of the estimates on the norm of the hyperinterpolation operators on the sphere, obtained by Sloan, Reimer, and others  (cf. \cite[Section~3.2]{sloansharma} for a review). We are also aware of a work \cite{damelin}  by Damelin and Levesley on similar questions on hyperinterpolation in the context of projective spaces.
   
The second application illustrates the use of \eref{mzgen1} when $p=\infty$. Suppose that the sampling nodes $\C_L$ are chosen randomly. Using probabilistic estimates, one is sometimes able to estimate the probability that the quantity $\sum_{y\in\C_L} |w_y| |\Phi_L(x,y)|$ exceeds a given threshold for any fixed $x\in\XX$. The inequality \eref{mzgen1} with $p=\infty$ then enables one to estimate the probability that the operator norm of $T_L^D$ is bounded. An example of this argument can be found in \cite{quadconst}. 

The third application concerns least square approximation. If one wishes to obtain $Q\in\Pi_L$ so as to minimize $\sum_{y\in\C_L} |w_y||f(y)-Q(y)|^2$, then one has to solve a system of linear equations with the matrix (the Gram matrix) $G$ whose $(j,k)$-th entry $G_{j,k}$ is given by $\sum_{y\in\C_L}|w_y|\phi_j(y)\phi_k(y)$. In view of the Rayleigh--Ritz theorem \cite[Theorem~4.2.2, p.~176]{horn}, the lowest  and highest eigenvalues of this matrix are given respectively by the infimum and supremum of the quotients
$$
\frac{\sum_{j,k}a_ja_kG_{j,k}}{\sum_{j}a_j^2}
$$
over all $a_j\in\RR$, $j=1,\cdots, \mbox{dim}(\Pi_L)$. Let $P=\sum_{j} a_j\phi_j$. Then the denominator expression above is equal to $\int_\XX |P(x)|^2d\mu(x)$. It is easy to check that the numerator expression is equal to $\sum_{y\in\C_L}|w_y||P(y)|^2$. Thus, if the weights $|w_y|$'s are chosen so as to  satisfy \eref{mzgen1} with $p=2$, then the lowest (respectively, the highest) eigenvalue is estimated from below (respectively, from above) by $c_1$ (respectively, $c_2$). In particular, the closer the ratio $c_2/c_1$ is to $1$, the better conditioned is the matrix $G$.

Finally, we have demonstrated in a number of papers starting with \cite{mnw1} that the inequalities \eref{mzgen1} with $p=1$ and $c_1,c_2$ sufficiently close to $1$ lead to the existence of positive quadrature formulas exact for integration of elements of $\Pi_{cL}$ for some constant $c$. 

Many modern applications require an analysis of huge, unstructured, high dimensional data sets, which are not dense on any cube, unlike classical approximation theory scenarios. Coifman and his collaborators have recently introduced diffusion geometry techniques for this purpose; see \cite{achaspissue} for an introduction. The idea is to assume that the data lies on an unknown low dimensional manifold. Lafon \cite{lafon} has shown that certain positive definite matrices constructed from the mutual distances among the data points with some tuning parameters converge to the ``heat kernel'' on the manifold. In theory, the heat kernel has a  formal representation of the form
\be\label{heatkerndefint}
K_t(x,y)=\sum_{k=0}^\infty \exp(-\ell_k^2t)\phi_k(x)\overline{\phi_k(y)},
\ee
where $\ell_k\uparrow \infty$ and the eigenfunctions $\{\phi_k\}$'s are orthonormal with respect to a suitable measure. 
This heat kernel has been used by Coifman and Maggioni \cite{mauro1} to define a metric on the data manifold, as well as to construct a multiresolution analysis on the manifold. In a more recent work \cite{jms}, Jones, Maggioni, and Schul have demonstrated that the heat kernel can be used to construct a local coordinate chart on the unknown manifold. 

Thus, even though the manifold is unknown, it is reasonable to assume for theoretical investigations that one knows a semi--group of positive definite kernels on the manifold, or equivalently (from a theoretical point of view) that one knows its infinitesimal generator, called the Laplacian on the unknown manifold.  
  In \cite{mauropap}, we started to develop a detailed theory of function approximation based on the eigenfunctions of the heat kernel. Our assumptions in \cite{mauropap} were formulated in terms of the behavior of the sums of the form $\sum_{\ell_k\le L}|\phi_k(x)|^2$ and the so called finite speed of wave propagation. However, since the actual manifold is unknown, and the heat kernel is the only easily computable quantity, we find it important in theoretical considerations to formulate the assumptions behind various theorems in terms of the heat kernel as far as possible. 
 In \cite{frankbern}, the properties of a summability kernel which plays a critical role in this theory were formulated purely in terms of the heat kernel, and generalized to obtain Marcinkiewicz--Zygmund inequalities, Bernstein inequalities, and the existence of quadrature formulas. 

This paper is devoted to a more detailed study of Marcinkiewicz--Zygmund inequalities in the case when $\XX$ is a smooth manifold, where the ``heat kernel'' based on the orthonormal system $\{\phi_j\}$ satisfies certain properties. In particular, our theory is valid when these are the eigenfunctions of the Laplace--Beltrami operator on the manifold, as well as in the case when they are eigenfunctions of certain weighted Laplace--Beltrami operators and a large class of other second order elliptic operators. We will establish conditions under which inequalities of the form \eref{mzgen1} hold for all $p$, $1\le p\le\infty$ with $c_1=1-\eta$ and $c_2=1+\eta$ for a prescribed $\eta$. Such inequalities were proved in \cite{frankbern} in the case when $p=1,\infty$, and applied to obtain the existence of quadrature formulas. However, a straightforward application of the classical Riesz--Thorin interpolation theorem is not sufficient to prove these inequalities in such a sharp form for $1<p<\infty$. We will prove an alternative form of this interpolation theorem, which appears to be new. We will also give intrinsic characterizations of the systems $\{(\C_L,W_L)\}$ without reference to the system $\{\phi_j\}$ which are equivalent to the inequalities of the form \eref{mzgen1} (without the requirement that $c_1, c_2$ should be arbitrarily close to $1$). For example, we will show that the upper inequality in \eref{mzgen1} holds if and only if for any $x\in\XX$ and any geodesic ball $B(x,r)$ of radius $r>0$ centered at $x$, 
$$
\sum_{y\in B(x,r)\cap \C_L} |w_y| \le c\mu(B(x,r+1/L)),
$$
for some constant $c>0$. Similar results will be proved for the lower inequality. Our aim is to provide such results for a very general setting, in particular including the case when the middle term in \eref{mzgen1} involves weighted averages of elements of $\Pi_L$ on balls rather than their values at finitely many points. In this paper, the heat kernel will play a somewhat indirect role. We will work with the very general setting of an arbitrary orthonormal system $\{\phi_k\}$ and sequence $\ell_k\uparrow\infty$. We will be using mainly the results in \cite{frankbern}. In turn, these are proved under the assumptions formulated in terms of the heat kernel defined formally by \eref{heatkerndefint}.
 
In Section~\ref{manifoldsect}, we will review a few basic facts regarding Riemannian manifolds in general, which will be needed in the rest of the paper. In Section~\ref{assumesect}, we discuss the various assumptions on the manifold, the measure, and the systems $\{\phi_j\}$, $\{\ell_j\}$ via the heat kernel. In Section~\ref{mzmeassect}, we introduce the abstract notions which enable us to generalize the theory from point evaluations to other measures. The main results are stated in Section~\ref{mainsect}, and proved in Section~\ref{proofsect}. These proofs require us to develop and review certain preparatory material, which is presented in  Section~\ref{prepsect}.

\bhag{Riemannian manifolds}\label{manifoldsect}
The purpose of this section is to review some facts and terminology regarding Riemannian manifolds. We will avoid very technical details, which can be found in such standard texts as \cite{boothby, docormo1, docormo2, oneil}. The material in this section is based mostly on \cite{docormo2}, and is essentially the same as the appendix to our paper \cite{frankbern}.

Let $q\ge 1$ be an integer. A \emph{differentiable manifold} of dimension $q$ is a set $\XX$ and a family of injective mappings $\x_\a :U_\a\subset \RR^q \to \XX$ of open sets $U_\a$ into $\XX$ such that (i) $\cup_{\a} \x_\a(U_\a)=\XX$, (ii) for any pair $\a,\beta$, with $\x_\a(U_\a)\cap \x_\beta(U_\beta)=W$ being nonempty, the sets $\x_\a^{-1}(W)$ and $\x_\beta^{-1}(W)$ are open subsets of $\RR^q$, and the mapping $\x_\beta^{-1}\circ \x_\a$ is (infinitely) differentiable on $\x_\a^{-1}(W)$, (iii) the family (\emph{atlas}) $\mathcal{A}_{\mathbb X}=\{(U_\a,\x_\a)\}$ is  maximal relative to the conditions (i) and (ii). The pair $(U_\a,\x_\a)$ (respectively, $\x_\a$) with $x\in \x_\a(U_\a)$ is called  a parametrization or coordinate chart (respectively, a system of coordinates) of $\XX$ around $x$, and $\x_\a(U_\a)$ is called a coordinate neighborhood of $x$. In the sequel, the term differentiable will mean infinitely many times differentiable. We assume also that $\XX$ is Hausdorff and has a countable basis as a topological space.

Intuitively, one thinks of a differentiable manifold as a surface in an ambient Euclidean space. The abstract definition above is intended to overcome the technical need for the ambient space. For all applications of our theory that we can imagine, and in particular, for an intuitive comprehension of our paper, there is no loss in thinking of a manifold as a surface. Moreover, a theorem of Whitney \cite[p.~30]{docormo2} provides a further justification of such a viewpoint.

Let $\mathbb{X},\mathbb{Y}$ be two differentiable manifolds. A mapping $f:\mathbb{X}\to\mathbb{Y}$ is called differentiable on an open set $W\subseteq\mathbb{X}$ if there exist for every $x\in W$ coordinate charts $(U,\x)\in\mathcal{A}_{\mathbb X}, (V,\y)\in\mathcal{A}_{\mathbb Y}$ with $x\in U, f(U)\subseteq V$ such that $\y^{-1}\circ f\circ \x$ is a $C^\infty$ function. In particular, a \emph{curve} in $\XX$ is a differentiable mapping from an interval in $\RR$ to $\XX$. The restriction of a curve $\gamma$ to a compact subinterval $[a,b]$ of $I$ is called a curve segment, joining $\gamma(a)$ to $\gamma(b)$. We may define a piecewise differentiable curve on a manifold $\XX$ in an obvious manner.

If $x\in \XX$, $\e>0$, and $\gamma :(-\e,\e)\to \XX$ is a curve with $x=\gamma(0)$, then the \emph{tangent vector} to $\gamma$ at $\gamma(t_0)$ is defined to be the functional $\gamma'(t_0)$ acting on the class of all differentiable $f:\XX\to \RR$ by
$$
\gamma'(t_0)f=\frac{d(f\circ\gamma)}{dt}\big\vert_{t=t_0}.
$$
The set of all such functionals $\gamma'(0)$ defines a vector space, called the \emph{tangent space} of $\XX$ at $x$, denoted by $T_x\XX$. Let $(U,\x)$ be a coordinate chart such that $0\in U$ and $x=\x(0)$, and for $j=1,\cdots, q$, $\partial_j(x)$ be the tangent vector at $x$ to the coordinate curve $x_j\to (0,\cdots,x_j,0,\cdots,0)$. Then $\{\partial_j(x)\}$ is a basis for $T_x\XX$. In particular, the dimension of $T_x\XX$ is $q$. The set $\{(x, v): x\in \XX, v\in T_x\XX\}$ is called the \emph{tangent bundle} of $\XX$, and can be endowed with the structure of a differentiable manifold of dimension $2q$. A \emph{vector field} $F$ on $\XX$ is a mapping that assigns to each $x\in\XX$ a vector $F(x)\in T_x\XX$ such that for every differentiable function $g$ on $\XX$, the mapping $x\mapsto F(x)g$ is differentiable.  If $G$ is another vector field, we may apply $G(x)$ to this mapping, obtaining thereby a second order vector field $G\circ F$. A derivative of higher order can be defined similarly.

A \emph{Riemannian metric} on a differentiable manifold $\XX$ is given by a scalar product $\langle \circ,\circ\rangle_x$ on each $T_x\XX$ which depends smoothly on the base point $x$, i.e. the function $\XX\to\RR$, $x\mapsto \langle X(x),Y(x)\rangle_x$ is $C^\infty(\XX)$. A manifold with a given Riemannian metric is called a \emph{Riemannian manifold}. Let $g_{i,j}=\langle \partial_{i}(x),\partial_{j}(x)\rangle_x$ and denote by $g$ the matrix $(g_{i,j})$. The entries of $g^{-1}$ are denoted be $g^{i,j}$.   The Riemannian metric on $\XX$ allows one to define a notion of length of a curve segment as well as the volume element (Riemannian measure) on $\XX$. First, if $F$ is a vector field on $\XX$, we may define $\tn F\tn_x := \langle F(x), F(x)\rangle_x$.
The length of a differentiable curve $\gamma :[0,1]\to \XX$ is defined as $\mathcal{L}(\gamma)=\int_0^1\tn \gamma'(t)\tn_{\gamma(t)} dt$. A differentiable curve $\gamma:[0,1]\to \XX$, such that the length of $\gamma$ does not exceed that of any other piecewise differentiable curve joining $\gamma(0)$ to $\gamma(1)$,  is called a \emph{geodesic} (\cite[Proposition~3.6, Corollary~3.9]{docormo2}). 

The gradient of a function $f\in C^\infty(\XX)$ is a vector field defined by 
$$
\nabla f =  \sum_{j=1}^q\sum_{i=1}^q g^{i,j}\partial_if\, \partial_{j}.
$$
For the gradient field we have 
\be\label{gradarbitrary}
\langle(\nabla f)_x,F(x)\rangle_x=(Ff)(x)
\ee
 for every vector field $F$. 
The divergence of a vector field $F=\sum_{j=1}^q F_j\partial_j$ is defined by
$$
\mbox{div }(F) =\frac{1}{\sqrt{\mbox{det}(g)}}\sum_{j=1}^q \partial_j(\sqrt{\mbox{det}(g)}F_j).
$$
The Laplace--Beltrami operator $\Delta^*f(x)$ is defined as the differential operator given by 
$$
\Delta^*f=-\mbox{div }(\nabla f)=\frac{-1}{\sqrt{\mbox{det}(g)}}\sum_{j=1}^q\sum_{i=1}^q\partial_j\left(\sqrt{\mbox{det}(g)}\, g^{i,j}\partial_i f\right). 
$$
The operator $\Delta^*$ is an elliptic operator. Therefore, in the case when $\XX$ is a compact connected manifold, the existence of a discrete spectrum and system of orthonormal eigenfunctions follows from the general theory of partial differential equations \cite[Chapter~5.1]{taylorbk}. 

We will have no further occasion to refer to the dimension of the manifold, and therefore, will use the symbol $q$ with different meanings in the rest of this paper.

\bhag{Assumptions}\label{assumesect}
Let $\XX$ be a compact, connected, Riemannian manifold (without boundary), $\rho$ be the geodesic distance on $\XX$, $\mu$ be a probability measure on $\XX$, and $\{\phi_k\}$ be an orthonormal system of continuous functions, $\phi_0(x)=1$ for all $x\in\XX$, and $\{\ell_k\}_{k=0}^\infty$ be an nondecreasing sequence of real numbers with $\ell_0=1$, $\ell_k\uparrow \infty$ as $k\to\infty$. For $L\ge 0$, the space $\span\{\phi_j : \ell_j\le L\}$ will be denoted by $\Pi_L$, and its members will be called \emph{diffusion polynomials} of degree at most $L$. We will also write $\Pi_\infty=\cup_{L\ge 0} \Pi_L$. For $x\in \XX$, $r>0$, let
$$
B(x,r):=\{y\in\XX\ :\ \rho(x,y)\le r\}, \ \Delta(x,r):=\XX\setminus B(x,r).
$$
\textbf{Assumption~1:}\\
\emph{
We assume that there exist  constants $\kappa_1, \a>0$ such that
\be\label{ballmeasurecond}
\mu(B(x,r))\le \kappa_1\, r^\a, \qquad x\in\XX, \ r>0.
\ee
}

Next, we discuss the notion of the heat kernel, and the assumptions on the same. 
For $t>0$, the \emph{heat kernel} is defined formally by
\be\label{heatkerndef}
K_t( x,y) =\sum_{k=0}^\infty \exp(-\ell_k^2t)\phi_k(x)\overline{\phi_k(y)}.
\ee
\textbf{Assumption 2:}\\
\emph{
We assume the existence of constants $\kappa_2,\kappa_3>0$ such that
\be\label{singlegaussbd}
|K_t(x,y)|\le \kappa_2t^{-\a/2}\exp(-\kappa_3\rho(x,y)^2/t), \qquad x,y\in\XX,\ t\in (0,1],
\ee
and
\be\label{singlegaussgradbd}
\tn\nabla_y K_t(x,y)\tn_x\le \kappa_2t^{-\a/2-1}\exp(-\kappa_3\rho(x,y)^2/t), \qquad x,y\in\XX,\ t\in (0,1],
\ee
where $\nabla_y$ indicates that the gradient is taken with respect to $y$.
We assume further that for some constant $\kappa_4>0$,
\be\label{lowergaussbd}
K_t(x,x)\ge \kappa_4 t^{-\a/2}, \qquad x\in\XX, \ t\in (0,1].
\ee
}
Kordyukov \cite{kordyukov91} has proved that each of our assumptions above hold for the heat kernels, when for each $k=0,1,\cdots$, $\phi_k$ is the eigenfunction  of a  second order elliptic operator corresponding to the eigenvalue $\ell_k^2$. The elliptic operator in question is assumed to satisfy some very general conditions, which are satisfied by the Laplace--Beltrami operator on a Riemannian manifold with ``bounded geometry'' (see \cite{kordyukov91} for definitions). Estimates on the heat kernel and its gradients are well understood in many other cases, including higher order partial differential operators on manifolds \cite{davies97, grigoryan95,  chengliyau81, dungey06}, with many other references given in \cite{grigoryan99}.

We will need  one more assumption, which we include here for the sake of organizational clarity, even though it requires some notations introduced in \eref{lpnormdef} and the paragraph which follows \eref{lpnormdef}.  A system $\{\psi_k\}\subset L^2(\mu)$  will be called a \emph{Bessel system} if there exists a dense subset ${\cal D}={\cal D}(\{\psi_k\})$ of  $C(\XX)$ (with respect to the norm of this space), such that (i) for any $\e>0$, ball of the form $B(x,r)$, and $f\in C(\XX) $ supported on $B(x,r)$, there exists $g\in {\cal D}$ such that the support of $g$ is contained in $B(x,2r)$ and 
$$
\int_\XX|f(x)-g(x)|d\mu(x)\le \e,
$$
 and (ii)
\be\label{besselineq}
\sum_{k=0}^\infty |\ip{f}{\psi_k}|^2 \le {\cal N}(f)<\infty, \qquad f\in {\cal D},
\ee
where ${\cal N}(f)$ is a positive number dependent on $f$, $\ip{\circ}{\circ}$ and $\{\psi_k\}$. 
Obviously, any orthonormal system on $\XX$ is a Bessel system with ${\cal D}= C(\XX)$. Another interesting example is the following. Let $\XX$ be a Riemannian manifold, $\mu$ be the Riemannian volume measure on $\XX$, $\{\phi_k\}$ be the eigenfunctions of the Laplace--Beltrami operator on $\XX$, $F$ be a  vector field on $\XX$, and $\psi_k=F\phi_k$, $k=0,1,\cdots$. For the space ${\cal D}$ we choose the class of all compactly supported, infinitely differentiable functions on $\XX$. Using  the Green's formula \cite[p.~383]{lee}, we obtain for any $f\in{\cal D}$, and $k=0,1,\cdots$,
\be\label{Topindentity}
\int_\XX \psi_k(x)f(x)d\mu(x)=-\int_\XX \phi_k(x)\mbox{ div }(fF)(x)d\mu(x).
\ee
Thus, $\{\psi_k\}$ is a Bessel family with ${\cal N}(f)=\|\mbox{ div }(fF)\|_{\mu;2}$.
A similar fact obtains also when one considers $wd\mu$ instead of $d\mu$ for some smooth positive valued function $w$.

\noindent\textbf{Assumption 3:}\\
\emph{
 We assume that for any  vector field $F$, the system $\{F\phi_k\}$ is a Bessel system.}

\noindent\textbf{Constant convention:}\\
\textit{In the sequel, the symbols $c, c_1,\cdots$ will denote positive constants depending only on $\XX$, $\rho$, $\mu$, $\kappa_1,\cdots,\kappa_4$, and other similar fixed quantities, but not on the systems $\{\phi_k\}$, $\{\ell_k\}$, nor any other variables not explicitly indicated. Their values may be different at different occurences, even within a single formula. The notation $A\sim B$ will mean $c_1A\le B\le c_2A$.
}

We note some consequences for our assumptions. We have proved \cite[Proposition~4.1]{frankbern}, \cite[Lemma~5.2]{eignet1} that the conditions \eref{singlegaussbd} with $x=y$ and \eref{lowergaussbd} are equivalent to
\be\label{christbd}
\sum_{\ell_j\le L}|\phi_j(x)|^2 \sim L^\a.
\ee
We note that the conditions $\phi_0(x)\equiv 1$, and $\ell_0=1$ imply that
\be\label{heatintegral}
\int_{\XX}K_t(x,y)d\mu(y)=1, \qquad x\in\XX, \ t\in (0,1].
\ee
In \cite{grigoryanheatmetric}, Grigor'yan has proved that \eref{ballmeasurecond}, \eref{heatintegral}, and  \eref{singlegaussbd} together imply that
\be\label{ballmeasurelowbd}
\mu(B(x,r))\ge cr^\alpha, \qquad 0<r\le 1, \ x\in\XX.
\ee
Using \eref{ballmeasurecond} and \eref{ballmeasurelowbd}, we  obtain that $\mu$ satisfies the homogeneity condition
\be\label{doublingcond}
\mu(B(x,R))\le c(R/r)^\a\mu(B(x,r)), \qquad x\in\XX,\ 0<r\le 1, R>0.
\ee


\bhag{MZ measures}\label{mzmeassect}
In this section, we wish to express the ideas in the introduction in a more abstract and formal manner.
First, it is cumbersome to write a sum of the form $\sum_{y\in\C}w_y f(y)$. To write such a sum, we need to introduce the set $\C$ and the weights $W=\{w_y\}$. The precise choice of these objects plays no role in our theory. Moreover, it makes it difficult to prescribe the dependence of various constants on the set $\C$ and the weights $W$. For these reasons we prefer to use the Lebesgue--Stieltjes integral notation $\int_\XX f(y)d\nu(y)$ to denote this sum, where $\nu$ is the measure that associates the mass $w_y$ with the point $y$, $y\in\C$; i.e., for $B\subset \XX$,
\be\label{countingmeasdef}
\nu(B)=\sum_{y\in B}w_y.
\ee
The notation has an additional advantage. When $f$ is not continuous, but in some $L^p(\mu)$, then $f$ cannot be defined everywhere. It is customary in such cases to consider averages of $f$ on small balls around the points in $\C$. A weighted sum of these averages can be written in the form $\int_\XX f(y)d\nu(y)$ as well, with a suitable choice of the measure $\nu$. Further, some applications require the consideration of a weighted manifold (cf. \cite{wtheatkernel}). Rather than dealing with the eigenfunctions for a weighted analogue of the Laplace--Beltrami operator, one may wish to work with the eigenfunctions for the unweighted case. The MZ inequalities with the measure $d\nu =wd\mu$, where $w$ is the weight function in question are expected to be useful in such situations. The class of all signed Borel measures on $\XX$ is a vector space, which will be denoted by ${\cal M}$. 

If $\nu\in {\cal M}$, its total variation measure is defined for Borel measurable subsets $B\subset\XX$ by 
$$
|\nu|(B) :=\sup \sum_{i=1}^\infty |\nu(B_i)|,
$$
where  the supremum is taken over all countable partitions $\{B_i\}$
of $B$. For any signed measure $\nu$, $|\nu|(\XX)$ is always a finite number. In the case when $\nu$ is the measure that associates the mass $w_y$ with each $y$, $y\in \C$,
one can easily deduce that $|\nu|(\XX)=\sum_{y\in \C}|w_y|$.
In the case when $d\nu=wd\mu$, one has $d|\nu|=|w|d\mu$. This includes the case when $\nu$ is a weighted sum of averages on balls. If $\nu\in {\cal M}$, the \emph{support} of $\nu$ is defined by 
$$
\supp(\nu):=\{x\in\XX : |\nu|(B(x,r))>0 \mbox{ for every } r>0\}.
$$
In view of \eref{ballmeasurelowbd}, $\supp(\mu)=\XX$.

Let $\nu$ be a signed measure on $\XX$. If $B\subseteq \XX$ is $\nu$-measurable, and $f : B\to \CC$ is a $\nu$-measurable function, we will write
\be\label{lpnormdef}
\|f\|_{\nu;p,B}:=\left\{\begin{array}{ll}
\disp \left\{\int_{B}|f(x)|^pd|\nu|(x)\right\}^{1/p}, & \mbox{ if $0< p<\infty$,}\\
|\nu|-\esssup_{x\in B}|f(x)|, &\mbox{ if $p=\infty$.}
\end{array}\right.
\ee
The class of all $f$ with $\|f\|_{\nu;p,B}<\infty$ will be denoted by $L^p(\nu;B)$, with the usual convention of considering two functions to be equal if they are equal $|\nu|$--almost everywhere. If $B=\XX$, we will omit its mention from the notations. The expressions $\|f\|_{\nu;p,B}$ are not norms if $p<1$, but we prefer to continue using the same notation. The inner product of $L^2(\mu)$ will be denoted by $\ip{\circ}{\circ}$. The $L^p(\mu)$--closure of $\Pi_\infty$ will be denoted by $X^p(\mu)$. The class of all uniformly continuous and bounded functions on $B$, equipped with the uniform norm will be denoted by $C(B)$. If $1< p<\infty$, the conjugate index $p'$ is defined by $p':=p/(p-1)$. We define $1'=\infty$ and $\infty'=1$.

Thus, if $\C=\C_L$ is as in the introduction, and $\nu$ is the corresponding measure as defined in \eref{countingmeasdef}, then the inequalities \eref{mzgen1} can be expressed in the concise form
\be\label{mzabstract}
c_1\|P\|_{\mu;p}\le \|P\|_{\nu;p}\le c_2\|P\|_{\mu;p}, \qquad P\in\Pi_L.
\ee

Since $\Pi_L$ is a finite dimensional space, such inequalities are always valid, but with the constants possibly depending on $L$. We are mostly interested in investigating the conditions under which these are independent of $L$, but wish to note another example where the constants depend polynomially on $L$. Since this provides an important rationale for considering general measures, apart from the wish to include averages over balls, we discuss this example in some detail. First, we note a property of diffusion polynomials, known as \emph{Nikolskii inequalities}.

\begin{prop}\label{nikolskiiprop}
Let $L>0$, $0<p<r\le \infty$, $P\in\Pi_L$. Then
\be\label{nikolskii}
\|P\|_{\mu;r} \le cL^{\a(1/p-1/r)}\|P\|_{\mu;p}.
\ee
\end{prop}
\begin{Proof}\ 
This proposition was proved in the case $p\ge 1$ in \cite[Lemma~5.5]{eignet1}. Let $0<p<1$, $P\in\Pi_L$. Then using the proved inequality \eref{nikolskii} with $r=\infty$ and $1$ in place of $p$, we obtain
$$
\|P\|_{\mu;\infty}\le cL^\a\int_\XX |P(x)|d\mu(x)= cL^\a\int_\XX |P(x)|^{1-p}|P(x)|^pd\mu(x)\le cL^\a\|P\|_{\mu;\infty}^{1-p}\|P\|_{\mu;p}^p.
$$
This leads to \eref{nikolskii} in the case when $r=\infty$ and $0<p<1$ as well. If $0<r<\infty$, then
$$
\|P\|_{\mu;r}^r =\int_\XX |P(x)|^{r-p}|P(x)|^pd\mu(x)\le \|P\|_{\mu;\infty}^{r-p}\|P\|_{\mu;p}^p\le cL^{\a(r-p)/p}\|P\|_{\mu;p}^{r-p}\|P\|_{\mu;p}^p= cL^{\a r (1/p-1/r)}\|P\|_{\mu;p}^r.
$$
This implies \eref{nikolskii} for all $r, p$, $0<p<r<\infty$.
\end{Proof}

\begin{example}\label{wtfnmzex}
Let $w \ge 0$ $\mu$--almost everywhere on $\XX$, $1<r<\infty$, and $w\in L^r(\mu)$. We define $d\nu=wd\mu$. Let $L>0$, $P\in\Pi_L$, and $0<p<\infty$.  Using H\"older inequality followed by \eref{nikolskii} with $pr'>p$ in place of $r$, we obtain
\begin{eqnarray*}
\|P\|_{\nu;p}^p &:=& \int_\XX |P(x)|^pw(x)d\mu \le \|w\|_{\mu;r}\disp\left\{\int_\XX |P(x)|^{pr'}d\mu(x)\right\}^{1/r'}\\
&=&\|w\|_{\mu;r}\|P\|_{\mu;pr'}^p \le cL^{\a p(1/p-1/(pr'))}\|w\|_{\mu;r}\|P\|_{\mu;p}^p= cL^{\a/r}\|w\|_{\mu;r}\|P\|_{\mu;p}^p.
\end{eqnarray*}
Thus, if $1<r<\infty$ and $w\in L^r(\mu)$, then
\be\label{ex1eqn1}
\|P\|_{\nu;p}\le c_1L^{\a/(pr)}\|w\|_{\mu;r}^{1/p}\|P\|_{\mu;p}.
\ee

To obtain an inequality in the reverse direction, let $1<q<\infty$, and $w^{-1}\in L^{q-1}(\mu)$. Using the Nikolskii inequality \eref{nikolskii} first with $p/q'$ in place of $p$ and $p$ in place of $r$, followed by H\"older inequality, we obtain
\begin{eqnarray*}
\|P\|_{\mu;p}^p&\le& c_2L^{\a(q'-1)}\|P\|_{\mu;p/q'}^p =c_2L^{\a/(q-1)}\disp\left\{\int_\XX |P(x)|^{p/q'}w^{1/q'}(x)w^{-1/q'}(x)d\mu(x)\right\}^{q'}\\
&\le& c_2L^{\a/(q-1)}\disp\left\{\int_\XX w^{-q/q'}(x)d\mu(x)\right\}^{q'/q}\left\{\int_\XX |P(x)|^pw(x)d\mu(x)\right\}\\
&=&c_2L^{\a/(q-1)}\|w^{-1}\|_{\mu;q-1}\|P\|_{\nu;p}^p.
\end{eqnarray*}
Thus, if $1<q<\infty$ and $w^{-1}\in L^{q-1}(\mu)$, then
\be\label{ex1eqn2}
\|P\|_{\mu;p}\le c_3L^{\a/(pq-p)}\|w^{-1}\|_{\mu;q-1}^{1/p}\|P\|_{\nu;p}.
\ee
In particular, if $w\in L^r(\mu)$ and $w^{-1}\in L^r(\mu)$, then
\be\label{ex1eqn3}
\|P\|_{\nu;p}\le c_1L^{\a/(pr)}\|w\|_{\mu;r}^{1/p}\|P\|_{\mu;p}, \quad \|P\|_{\mu;p}\le c_3L^{\a/(pr)}\|w^{-1}\|_{\mu;r}^{1/p}\|P\|_{\nu;p}.
\ee
\end{example}

The measure $\nu$ will be called an \emph{MZ measure} if the constants $c_1,c_2$ appearing in \eref{mzabstract} are independent of $L$.

\bhag{Main theorems}\label{mainsect}
Let $\C\subset K\subset \XX$ be compact sets. We define the mesh norm $\delta(\C,K)$ of $\C$ with respect to $K$ and the minimal separation of $\C$ by
\be\label{meshnormdef}
\delta(\C,K) =\sup_{x\in K}\rho(x,\C),\qquad q(\C)=\min_{x,y\in\C,\ x\not=y}\rho(x,y).
\ee
To keep the notation simple, we will write $\delta(\C):=\delta(\C,\XX)$. Of course, the quantity $q(\C)$ is of interest only when $\C$ is a finite set. It is easy to see that $q(\mathcal{C})/2\leq\delta(\mathcal{C})$. 

Our first theorem states the MZ inequalities in a sharp form in an apparently special  case. We note that part (a) of the following theorem was proved (with minor differences) in \cite[Theorem~3.2]{frankbern} for the case $p=1$.

\begin{theorem}\label{mztheolp}
Let $\C=\{x_1,\cdots,x_M\}$ be a finite subset of $\XX$ satisfying 
\be\label{cuniformitygen}
\frac{1}{2}q(\C)\le \delta(\C)\le \kappa q(\C)
\ee
for some $\kappa\geq 1$.  
Let $1\le p<\infty$ and $A\ge 2$. In this theorem, all the constants may depend upon $\kappa$ and $A$.\\
{\rm (a)} There exist $c_1, c_2>0$ such that for every $\eta>0$, if $\delta(\C)\le c_1$,  $L\le c_2\eta(p\delta(\C))^{-1}$, and $P\in\Pi_L$, then
\be\label{mzineqstronglp}
 \sum_{k=1}^M\mu(B(x_k,\delta(\C)))\sup_{z, y\in B(x_k,A\delta(\C))}\left||P(z)|^p-|P(y)|^p\right|\le \eta\int_\XX |P(x)|^pd\mu(x),
\ee
{\rm (b)} Let $\C$ be as in part (a), and $\{Y_k\}_{k=1}^M$ be a partition of $\XX$ 
  such that $x_k\in Y_k\subseteq B(x_k,A\delta(\C))$ for each $k$, $1\le k\le M$. There exists $c_3>0$ such that for $L\le c_3\eta(p\delta(\C))^{-1}$, and $P\in\Pi_L$, we have
\be\label{discstronglpmzineq}
\left|\int_\XX |P(z)|^pd\mu(z)-\sum_{k=1}^M \mu(Y_k)|P(x_k)|^p\right| \le \eta\int_\XX |P(x)|^pd\mu(x).
\ee
{\rm (c)} There exists $c_4>0$ such that if $L \le c_4\eta\delta(\C)^{-1}$ then
\be\label{mzineqstronginf}
\left|\|P\|_{\mu;\infty}-\max_{1\le k\le M}|P(x_k)|\right| \le \eta\|P\|_{\mu;\infty}.
\ee
\end{theorem}

We note that a variant of Theorem~\ref{mztheolp} was stated in our paper \cite[Theorem~3.1]{mnw1} in the case when $\XX$ is the Euclidean unit sphere, $\phi_j$'s are spherical harmonics (so that $\Pi_L$ is the class of all spherical polynomials of degree at most $L$), and $\mu$ is the Riemannian volume measure on the sphere. The theorem is correct for $p=1,\infty$ as stated there, but the proof does not use the correct form of the Riesz--Thorin interpolation theorem which is needed for proving such inequalities.  Also, in the proof of \cite[Theorem~3.2]{frankbern}, we had constructed a partition $Y_k$. However, it was an error on our part to assume that $x_k\in B(x_k,q(\C)/2)\subset Y_k$. Both of these errors are corrected in Theorem~\ref{mztheolp} and the proof of Theorem~\ref{newstrongmztheo}.

Next, we wish to give an analogue of Theorem~\ref{mztheolp} where general measures are involved. The transition from the finitely supported measures to the general case is achieved via the following theorem.

\begin{theorem}\label{genpartitiontheo}
Let $\nu$ be a signed measure, $\delta(\supp(\nu))<d\le 1/81$. Then there exists a finite subset $\C=\{x_1,\cdots,x_M\}\subseteq \supp(\nu)$ with the property that
\be\label{subsetdelta}
q(\C)/2\le \delta(\C)\le 81d \le 162q(\C).
\ee
Moreover, there exists a partition $\{Y_k\}_{k=1}^M$ of $\XX$ and a finite subset $\tilde \C$ with $\C\subseteq \tilde\C\subseteq \supp(\nu)$ such that for  $k=1,\cdots, M$, $x_k\in Y_k\subseteq B(x_k, 81d)$, $\mu(Y_k)\sim d^\a$, and $|\nu|(Y_k)\ge c\min_{x\in \tilde \C}|\nu|(B(x,d/4))>0$. 
\end{theorem}

Theorem~\ref{genpartitiontheo} helps us to use Theorem~\ref{mztheolp} to arrive at the following statement, where general measures are involved.

\begin{theorem}\label{newstrongmztheo}
Let $\nu$ be a signed measure, $\delta(\supp(\nu))<d\le 1/81$, $\C$, $\{Y_k\}$ be as in Theorem~\ref{genpartitiontheo}.  \\
{\rm (a)}  Let $1\le p<\infty$. There exist $c_1, c_2>0$ such that for every $\eta>0$, if $d\le c_1$,  $L\le c_2\eta(pd)^{-1}$, and $P\in\Pi_L$, then
\be\label{newstrongmzineq}
\sum_{x\in\C} \left|\int_{Y_k}|P(y)|^pd\mu(y)-\frac{\mu(Y_k)}{|\nu|(Y_k)}\int_{Y_k}|P(z)|^pd|\nu|(z)\right| \le \eta\|P\|_{\mu;p}^p.
\ee
{\rm (b)} There exists $c_3>0$ such that if $L\le c_3\eta d^{-1}$ and $P\in\Pi_L$, then
\be\label{newsupmzineq}
\left|\|P\|_{\nu;\infty} -\|P\|_{\mu;\infty}\right| \le \eta\|P\|_{\mu;\infty}.
\ee
\end{theorem}

Since  Theorem~\ref{newstrongmztheo}(b) settles the question of MZ inequalities in the case $p=\infty$,  we will focus in the remainder of this paper on the case when $1\le p<\infty$. It is clear from Theorem~\ref{newstrongmztheo}(a) that the MZ inequalities for the measure $\nu$ will depend upon the relationship between $\nu(B(x,d))$ and $\mu(B(x,d))\sim d^\a$ for $x\in\XX$. Accordingly, we make the following definition.

\begin{definition}\label{regularitydef}
Let $\nu\in {\cal M}$, $d>0$.\\
{\rm (a)} We say that $\nu$ is \textrm{$ d$--regular} if
\be\label{regulardef}
\nu(B(x,d))\le cd^\a, \qquad x\in\XX.
\ee
The infimum of all constants $c$ which work in \eref{regulardef} will be denoted by $\tn\nu\tn_{R,d}$. \\
{\rm (b)} We say that $\nu$ is \textrm{$d$--dominant} if
\be\label{dominantdef}
\nu(B(x,d))\ge cd^\a, \qquad x\in\XX.
\ee
The supremum of all $c$ which work in \eref{dominantdef} will be denoted by $\tn\nu\tn_{D,d}^{-1}$. 
\end{definition}
We observe that \eref{regulardef} and \eref{dominantdef} are very similar to \eref{ballmeasurecond} and \eref{ballmeasurelowbd} respectively. However, in contrast to  \eref{ballmeasurecond} and \eref{ballmeasurelowbd}, \eref{regulardef} and \eref{dominantdef} are required to hold only for one value of $d$. Also, the function $\nu\to \tn\nu\tn_{R,d}$ is a norm on the space of all $d$--regular measures.

\begin{example}\label{discsetex}

 Let $\C$ be as in Theorem~\ref{mztheolp}, $\nu$ be the measure that associates the mass $\mu(B(x_k,\delta(\C)))$ with each $x_k$, $k=1,\cdots,M$. Let $x\in\XX$, $\tilde\C=B(x,2\delta(\C))\cap \C$. Then \eref{cuniformitygen} implies that the balls $B(y,\delta(\C)/(2\kappa))$, $y\in\tilde\C$ are mutually disjoint, and clearly, their union is a subset of $B(x,2\delta(\C)(1+1/(2\kappa)))$. So, in view of \eref{ballmeasurecond} and \eref{doublingcond}, we obtain
\bea\label{preregcond}
\nu(B(x,2\delta(\C))) &=&\sum_{y\in\tilde\C}\mu(B(y,\delta(\C)))\le c\sum_{y\in\tilde\C}\mu(B(y,\delta(\C)/(2\kappa))) =c\mu\left(\cup_{y\in\tilde\C}B(y,\delta(\C)/(2\kappa))\right)\nonumber\\
&\le& c\mu(B(x,2\delta(\C)(1+1/(2\kappa))))\le c_1\delta(\C)^\a.
\eea
In the reverse direction, the definition of the mesh norm implies that $B(x,\delta(\C))\subset \cup_{y\in\tilde\C}B(y,\delta(\C))$. Therefore, we obtain using \eref{ballmeasurelowbd} that
\be\label{predomcond}
\nu(B(x,2\delta(\C)))=\sum_{y\in\tilde\C}\mu(B(y,\delta(\C)))\ge \mu(B(x,\delta(\C)))\ge c\delta(\C)^\a.  
\ee
Thus, $\nu$ is $2\delta(\C)$--regular as well as $2\delta(\C)$--dominant.
\end{example}

\begin{example}\label{capavgex}
 Let $\C$ be as in Theorem~\ref{mztheolp}. For each $y\in\C$, let $q(\C)/4<r_y\le q(\C)/2$, $\chi_y$ be the characteristic (indicator) function of $B(y,r_y)$. Then using the same argument as above, it is easy to verify that the measure $d\nu =(\sum_{y\in\C}\chi_y)d\mu$ also satisfies  \eref{preregcond} and \eref{predomcond} (with different constants). Thus, this $\nu$ is also $c_1\delta(\C)$--regular and $c_2\delta(\C)$--dominant.
\end{example}

\begin{example}\label{wtfnregex}
Let $w\in L^r(\mu)$ for some $1<r<\infty$, $w\ge 0$ $\mu$--almost everywhere on $\XX$, and $d\nu=wd\mu$. Then $\supp(\nu)=\XX$. Let $x\in\XX$, $d\in (0,2/3)$, and in this example only, $M=\|w\|_{\mu;r}/d^{\a/r}$, $E=\{y\in B(x,d) : w(y)\ge M\}$. Then H\"older's inequality implies that
$$
M\mu(E) \le \int_E wd\mu \le \|w\|_{\mu;r}\mu(E)^{(r-1)/r},
$$
and hence, $\mu(E)\le (\|w\|_{\mu;r}/M)^r$. The second inequality above then yields
$$
\int_E wd\mu \le \|w\|_{\mu;r}^rM^{-(r-1)}.
$$
Therefore, our choice of $M$ implies that
\be\label{ex2eqn1}
\nu(B(x,d)) =\int_{B(x,d)\setminus E}w d\mu +\int_E wd\mu \le M\mu(B(x,d))+ \|w\|_{\mu;r}^rM^{-(r-1)} \le (\kappa_1+1)M d^\a,
\ee
where $\kappa_1$ is defined in \eref{ballmeasurecond}. 
To obtain an estimate analogous to \eref{dominantdef}, let $w^{-1}\in L^{q-1}(\mu)$ (i.e, $w^{-1}\in L^q(\nu)$) for some $q>1$. The same argument as above shows that for any $M_1>0$, $\tilde E=\{y\in B(x,d) : w(y)^{-1}\ge M_1\}$, we have
$$
\int_{\tilde E}w^{-1} d\nu =\mu(\tilde E) \le \|w^{-1}\|_{\nu;q}^qM_1^{-(q-1)}= \left(\frac{\|w^{-1}\|_{\mu;q-1}}{M_1}\right)^{q-1}.
$$
Hence, in view of \eref{ballmeasurelowbd},
$$
cd^\a \le \mu(B(x,d)) =\int_{B(x,d)\setminus \tilde E}w^{-1}d\nu + \int_{\tilde E}w^{-1}d\nu \le M_1\nu(B(x,d)) + \left(\frac{\|w^{-1}\|_{\mu;q-1}}{M_1}\right)^{q-1}.
$$
We now choose $M_1=\|w^{-1}\|_{\mu;q-1}(cd^\a/2)^{-1/(q-1)}$, and conclude that
\be\label{ex2eqn3}
\nu(B(x,d))\ge c_1M_1^{-1}d^\a.
\ee
\end{example}

\begin{theorem}\label{mzmeascharacttheo}
Let $L\ge 2$,  and $\nu\in {\cal M}$. In this theorem, all constants $c, c_1,\cdots$, may depend upon $p$.\\
{\rm (a)} If $\nu$ is $1/L$--regular, then $\|P\|_{\nu;p}\le c_1\tn\nu\tn_{R,1/L}^{1/p}\|P\|_{\mu;p}$ for all $P\in \Pi_L$ and $1\le p<\infty$. Conversely, if for some $A>0$ and $1\le p<\infty$, $\|P\|_{\nu;p}\le A^{1/p}\|P\|_{\mu;p}$ for all $P\in \Pi_L$, then $\nu$ is $1/L$--regular, and $\tn\nu\tn_{R,1/L}\le c_2A$. \\
{\rm (b)} There exists  constants $c, c_4$ such that if $L\ge c$ and $\nu$ is  $c_4/L$--dominant, then $\|P\|_{\mu;p}\le c_3\tn\nu\tn_{D,c_4/L}^{1/p}\|P\|_{\nu;p}$ for all $P\in \Pi_L$, and for all $p$, $1\le p<\infty$. Conversely, let $\nu$ be $1/L$--regular, and $S>\a$ be an integer. If for some $A_1>0$, and $1\le p<\infty$, $\|P\|_{\mu;p}\le A_1^{1/p}\|P\|_{\nu;p}$ for all $P\in \Pi_L$, then $\nu$ is $d=c_5(S)(\max(1,\tn\nu\tn_{R,1/L}A_1)^{1/(S-\a)})L^{-1}$--dominant, and $\tn\nu\tn_{D,d}\ge c_6(S)A_1$.
\end{theorem}

The term $d$--regular has been used with different meanings in our other papers. The following proposition reconciles the different definitions.

\begin{prop}\label{mzequivprop}
Let $d\in (0,1]$, $\nu\in {\cal M}$. \\
{\rm (a)} If $\nu$ is $d$--regular, then for each $r>0$ and $x\in\XX$,
\be\label{regreconcile}
|\nu|(B(x,r))\le c\tn\nu\tn_{R,d}\ \mu(B(x,r+d))\le  c_1\tn\nu\tn_{R,d}(r+d)^\a.
\ee
Conversely, if for some $A>0$, $|\nu|(B(x,r))\le A(r+d)^\a$ or each $r>0$ and $x\in\XX$, then $\nu$ is $d$--regular, and $\tn\nu\tn_{R,d}\le 2^\a A$.\\
{\rm (b)} For each $\gamma>1$, 
\be\label{regequiv}
\tn\nu\tn_{R,\gamma d}\le c_1(\gamma+1)^\a \tn\nu\tn_{R,d}\le c_1(\gamma+1)^\a\gamma^\a\tn\nu\tn_{R,\gamma d},
\ee
where $c_1$ is the constant appearing in \eref{regreconcile}.
\end{prop}

We end this section with an discussion about positive quadrature formulas. We will say that $\nu$ is a \emph{quadrature measure of order $L$} if 
$$
\int_\XX P(y)d\mu(y)=\int_\XX P(y)d\nu(y), \qquad P\in \Pi_L.
$$

First, we prove a very general existence theorem for such formulas.

\begin{theorem}\label{quadtheo}
There exist constants $c_1, c_2>0$ with the following property: If $\nu$ is a signed measure, $\delta(\supp(\nu))< d<c_1$ and $0<L< c_2d^{-1}$, then there exists a simple function $W :\supp(\nu)\to [0,\infty)$, satisfying 
\be\label{diffpolyquad}
\int_\XX P(y)d\mu(y) =\int_\XX P(y)W(y)d|\nu|(y), \qquad P\in\Pi_L.
\ee
If $\nu$ is $d$--regular, then $W(y)\ge c\tn\nu\tn_{R,d}^{-1}$, $y\in\XX$. 
\end{theorem}

We observe that if $\nu$ is supported on a finite subset of $\XX$, then this reduces to \cite[Theorem~3.1(b)]{frankbern}.  We find it remarkable that the only conditions on $\nu$ for \eref{diffpolyquad} to hold are on $\supp(\nu)$. 

In many cases of interest, for example, the Euclidean sphere, the rotation group $SO(3)$ and projective spaces, if $P\in \Pi_L$, then $P^2\in \Pi_{2L}$. In the appendix, we will show that a similar fact holds for eigenfunctions of a fairly large class of elliptic operators. In the very general situation considered in this paper, we make the following product assumption as in \cite{eignet1}. To formulate this assumption, we need one further notation. If $f\in L^p(\mu)$, and $m>0$, we denote
$$
\dist(p;f,\Pi_m):=\inf_{P\in\Pi_m}\|f-P\|_{\mu;p}.
$$

\noindent\textbf{Product assumption:}\\
\emph{
We assume that there exists a constant $A^*\ge 2$ with the following property: With
\be\label{gammanormdef}
\e_L:=\sup_{\ell_j,\ell_k\le L}\dist(\infty;\phi_j\phi_k,\Pi_{A^*L}), \qquad L>0,
\ee
 we have $L^c\e_L\to 0$ as $L\to\infty$ for every $c>0$. }\\

We have conjectured in \cite{eignet1} that this assumption holds for every analytic manifold $\XX$.

\begin{theorem}\label{posmztheo}
Let the product assumption hold. There exists a constant $c>0$ such that if $L\ge c$ and $\tau$ is a positive quadrature measure of order $2A^*L$, then  
\be\label{wdnumzineq}
\|P\|_{\tau;p}\sim \|P\|_{\mu;p}, \qquad P\in \Pi_L,
\ee
where the constants involved may depend upon $p$ but not on $\tau$ or $L$.
\end{theorem}

\bhag{Preparatory results}\label{prepsect}
In this section, we summarize some results which will be needed in the proofs of the theorems in Section~\ref{mainsect}. In Section~\ref{rieszsect}, we prove the Riesz--Thorin interpolation theorem in the form in which we need it.  In Section~\ref{locopsect}, we summarize some of the properties of a localized kernel and diffusion polynomials \cite{mauropap, frankbern}, and  extend these to the $L^p$ setting using the Riesz--Thorin interpolation theorem. In Section~\ref{kreinsect}, we prove, for the sake of completeness, a special case of Krein's extension theorem for positive functionals, following a hint in \cite[Exercise~(14.27), p.~200]{hewittbk}. This will be used in proving the existence of positive quadrature formulas in Theorem~\ref{quadtheo}(b). 

\subsection{Riesz--Thorin interpolation theorem}\label{rieszsect}
Let $\X$, $\Y$ be Banach spaces of functions defined on a measure space $(\Omega,\tau)$. We assume the existence of associated Banach spaces $\X'$, $\Y'$ such that
\be\label{duality}
\|f\|_\X=\sup\left\{\int_\Omega f(x)\overline{g(x)}d\tau(x) : \|g\|_{\X'}=1\right\}, \quad  \|f\|_\Y=\sup\left\{\int_\Omega f(x)\overline{g(x)}d\tau(x) : \|g\|_{\Y'}=1\right\}.
\ee
Let ${\bf W}=\{w_k\}_{k=1}^M\subset (0,\infty)$. For $1\le p\le \infty$,  and integer $M\ge 1$, we define for ${\bf a}=(a_1,\cdots,a_M)\in\CC^M$,
$$
\|{\bf a}\|_{{\bf W},\ell^p} =\left\{\begin{array}{ll}
\left(\sum_{k=1}^M w_k|a_k|^p\right)^{1/p}, & \mbox{ if $1\le p<\infty$,}\\
\max_{1\le k\le M}|a_k|, &\mbox{ if $p=\infty$}.
\end{array}\right.
$$
 It is elementary to check that
\be\label{littlelpduality}
\|{\bf a}\|_{{\bf W},\ell^p}=\sup\left\{\sum_{k=1}^M w_ka_k\overline{b_k}  : \|(b_1,\cdots,b_M)\|_{{\bf W},\ell^{p'}}=1\right\}.
\ee
 We define the space $\X_{{\bf W},p}$ to be the tensor product space $\otimes_{k=1}^M \X$ equipped with the norm
$$
\|{\bf f}\|_{\X, {\bf W}, p} :=\|(\|f_1\|_\X,\cdots,\|f_M\|_\X)\|_{{\bf W},\ell^p}, \qquad {\bf f}=(f_1,\cdots,f_M)\in \X_{{\bf W},p}.
$$
The space $\X'_{{\bf W},p}$ and the norm $\|\circ\|_{\X',{\bf W},p}$ are defined similarly.

In the statement of the Riesz--Thorin interpolation theorem, we need another measure space. Not to complicate our notations, we will use $(\XX,\mu)$ here. However, it should be understood that the only property we need is that this is a  measure space, with $\mu$ being a positive measure. In this subsection we are not assuming any properties and other assumptions, including the fact that $\XX$ is a manifold, and $\mu$ is a probability measure.

\begin{theorem}\label{tensorrieszthorintheo}
Let $1\le p_0\le   p_1\le \infty$, $1\le r_0\le r_1\le\infty$, $0<t<1$, ${\mathcal U}$ be an operator satisfying
\be\label{uopboundarycond}
\|{\mathcal U}f\|_{\X, {\bf W}, p_j}\le M_j\|f\|_{\mu;r_j}, \qquad f\in L^{r_j}(\mu), \ j=1,2,
\ee
and $1/p=(1-t)/p_0 +t/p_1$, $1/r=(1-t)/r_0 +t/r_1$. Then 
\be\label{uopbdlp}
\|{\mathcal U}f\|_{\X, {\bf W}, p}\le M_0^{1-t}M_1^t\|f\|_{\mu;r}, \qquad f\in L^r(\mu).
\ee
\end{theorem}
The proof of Theorem~\ref{tensorrieszthorintheo} mimicks that of the usual Riesz--Thorin theorem. We could not find a reference where this theorem is stated in the form in which we need it. Therefore, we include a proof, following that of the usual Riesz--Thorin theorem as given in \cite[Chapter~XII, Theorem~1.11]{zygmund}. The first step is the following lemma.

 \begin{lemma}\label{dualitylemma}
\be\label{tensorduality}
\|{\bf f}\|_{\X, {\bf W}, p}=\sup\left|\sum_{k=1}^Mw_k \overline{b_k}\int_\Omega f_k(x)\overline{g_k(x)}d\tau(x)\right|,
\ee
where the supremum is over all ${\bf b}=(b_1,\cdots,b_M)$ with $\|{\bf b}\|_{{\bf W},\ell^{p'}}=1$ and $g_1,\cdots,g_M\in \X'$ with $\|g_k\|_{\X'}=1$, $k=1,\cdots,M$.
\end{lemma}
\begin{Proof} 
In view of H\"older's inequality, it is clear that
\be\label{pf1eqn1}
\sup\sum_{k=1}^Mw_k \overline{b_k}\int_\Omega f_k(x)\overline{g_k(x)}d\tau(x)\le \|{\bf f}\|_{\X, {\bf W}, p},
\ee
where the supremum is over all ${\bf b}$ with $\|{\bf b}\|_{{\bf W},\ell^{p'}}=1$ and $g_1,\cdots,g_M\in \X'$ with  $\|g_k\|_{\X'}=1$, $k=1,\cdots,M$. If ${\bf f}=0$ then \eref{tensorduality} is obvious. Let ${\bf f}\not=0$, and $\e>0$. In view of \eref{duality} and \eref{littlelpduality}, there exist $g_k\in\X'$ and ${\bf b}\in [0,\infty)^M$ such that $\|g_k\|_{\X'}=1$, $k=1,\cdots,M$, $\|{\bf b}\|_{{\bf W},\ell^{p'}}=1$ and
$$
\sum_{k=1}^Mw_k \overline{b_k}\int_\Omega f_k(x)\overline{g_k(x)}d\tau(x)\ge (1-\e)\sum_{k=1}^Mw_k \overline{b_k}\|f_k\|_\X\ge (1-\e)^2\|{\bf f}\|_{\X,{\bf W},p}.
$$
\end{Proof}

Next, we recall the Phragm\'en--Lindel\"of maximum principle \cite[Chapter~XII, Theorem~1.3]{zygmund}.
\begin{prop}\label{threelineprop}
Supose that $f$ is continuous and bounded on the closed strip of the complex plane $0 \le \Re e\ z\le 1$, and analytic in the interior of this strip. If $|f(z)|\le M_0$,  $\Re e\  z=0$,  and $|f(z)|\le M_1$, $\Re e\  z=1$, then $|f(z)|\le M_0^{1-t}M_1^t$ for $\Re e\  z = t$.
\end{prop}

We are now ready to prove Theorem~\ref{tensorrieszthorintheo}.

\noindent
\textsc{Proof of Theorem~\ref{tensorrieszthorintheo}.} In this proof only, we will write $\a_j=1/r_j$, $\beta_j=1/p_j$, $j=1,2$,  $\a(z)=(1-z)\a_0+z\a_1$, $\beta(z)=(1-z)\beta_0+z\beta_1$, so that $\a(t)=1/r$, $\beta(t)=1/p$. If $r_0=r_1=\infty$, then $r=\infty$ as well, and \eref{uopbdlp} is a simple consequence of H\"older inequality. So, we assume that $r_0<\infty$, and hence, $r<\infty$. Since simple functions are dense in $L^r(\mu)$, it is enough to prove \eref{uopbdlp} when $f$ is a simple function; i.e., $f=\sum_{j=1}^N d_je^{iu_j}\chi_j$, where $N\ge 1$ is an integer, $d_j>0$, $u_j\in (-\pi,\pi]$, and $\chi_j$'s are the characteristic functions of pairwise disjoint sets. We define
$f_z=\sum_{j=1}^N d_j^{\alpha(z)r}e^{iu_j}\chi_j$. Next, let ${\bf b}=(|b_1|e^{iv_1},\cdots,|b_M|e^{iv_M})\in\CC^M$ be an arbitrary vector satisfying $\|{\bf b}\|_{{\bf W},\ell^{p'}}=1$, and $g_1,\cdots,g_M\in \X'$ be arbitrary functions satisfying $\|g_k\|_{\X'}=1$, $k=1,\cdots,M$. We define 
$$
G_{z,k}=|b_k|^{(1-\beta(z))p'}e^{-iv_k}, \qquad k=1,\cdots,M,
$$
where it is understood that $G_{z,k}=0$ if $b_k=0$. Finally, we define
$$
\Phi(z)=\sum_{k=1}^M w_kG_{z,k}\int_\Omega ({\mathcal U}f_z)_k(x)\overline{g_k(x)}d\tau(x)=\sum_{k=1}^M\sum_{j=1}^N w_kG_{z,k}d_j^{\alpha(z)p}e^{iu_j}\int_\Omega ({\mathcal U}\chi_j)(x)\overline{g_k(x)}d\tau(x).
$$
We note that $f_t=f$, $G_{t,k}=\overline{b_k}$, and therefore,
\be\label{pf2eqn2}
\Phi(t)=\sum_{k=1}^Mw_k\overline{b_k}\int_\Omega ({\mathcal U}f)_k(x)\overline{g_k(x)}d\tau(x).
\ee
Now, $\Phi$ is a finite linear combination of functions of the form $e^{az}$, and hence, is an entire function, bounded on the strip $0\le \Re e\  z \le 1$.  If $\Re e\  z=0$ then $|G_{z,k}|=|b_k|^{p'\Re e\  (1-\beta(z))}=|b_k|^{p'(1-\beta_0)}=|b_k|^{p'/p_0'}$.   Therefore, using H\"older's inequality, we obtain for $\Re e\  z=0$:
\bea\label{pf2eqn1}
|\Phi(z)|&\le& \|(G_{z,1},\cdots,G_{z,M})\|_{{\bf W},\ell^{p_0'}}\left\|\left(\int_\Omega ({\mathcal U}f_z)_1(x)\overline{g_1(x)}d\tau(x), \cdots, \int_\Omega ({\mathcal U}f_z)_M(x)\overline{g_M(x)}d\tau(x)\right)\right\|_{{\bf W},\ell^{p_0}}\nonumber\\
&\le&\|{\bf b}\|_{{\bf W},\ell^{p'}}^{p'/p_0'}\|(\|({\mathcal U}f_z)_1\|_\X,\cdots, \|({\mathcal U}f_z)_M\|_\X)\|_{{\bf W}, \ell^{p_0}}= \|{\mathcal U}f_z\|_{\X,{\bf W},p_0}\le M_0\|f_z\|_{\mu;r_0}.
\eea
For $\Re e\  z=0$, $|d_j^{\a(z)r}|=d_j^{r/r_0}$. Also, at any point $x\in\XX$, there is at most one $j$ such that $\chi_j(x)\not=0$. For this $j$, $f_z(x)= d_j^{\a(z)r}e^{iu_j}\chi_j(x)$, and $f(x)=d_je^{iu_j}\chi_j(x)$. So, for $\Re e\  z=0$, and any $x\in\XX$, $|f_z(x)|^{r_0} =\sum_{j=1}^N d_j^r\chi_j(x)=|f(x)|^r$. Thus, $\|f_z\|_{\mu;r_0}=\|f\|_{\mu;r}^{r/r_0}$. Hence, \eref{pf2eqn1} shows that $|\Phi(z)|\le M_0\|f\|_{\mu;r}^{r/r_0}$, $\Re e\  z=0$. Similarly, $|\Phi(z)|\le M_1\|f\|_{\mu;r}^{r/r_1}$, $\Re e\  z=1$. Proposition~\ref{threelineprop} then implies that $|\Phi(t)|\le M_0^{1-t}M_1^t\|f\|_{\mu;r}$; i.e., in view of \eref{pf2eqn2}, we have
$$
 \left|\sum_{k=1}^Mw_k\overline{b_k}\int_\Omega ({\mathcal U}f)_k(x)\overline{g_k(x)}d\tau(x)\right| \le M_0^{1-t}M_1^t\|f\|_{\mu;r}.
$$
Since ${\bf b}$ and the functions $g_k$ were arbitrary subject only to $\|{\bf b}\|_{{\bf W},p'}=1$, $\|g_k\|_{\X'}=1$, the estimate \eref{uopbdlp} follows from Lemma~\ref{dualitylemma}.
\qed

\subsection{Localized polynomial operators}\label{locopsect}
Let $h :\RR\to [0,1]$ be an even, $C^\infty$ function, nonincreasing on $[0,\infty)$ such that $h(t)=1$ if $|t|\le 1/2$ and $h(t)=0$ if $|t|\ge 1$. We will treat $h$ to be a fixed function, so that the dependence of different constants on the choice of $h$ will not be mentioned. We will write
\be\label{phikerndef}
\Phi_L(x,y):=\sum_{j=0}^\infty h(\ell_j/L)\phi_j(x){\phi_j(y)}.
\ee
For $f\in L^1(\mu)$, we define
$$
\hat f(k)=\int_\XX f(y){\phi_k(y)}d\mu(y), \qquad k=0,1,\cdots,
$$
and
\be\label{sigmaopdef}
\sigma_L(f,x):=\int_\XX f(y)\Phi_L(x,y)d\mu(y) =\sum_{\ell_j\le L}h(\ell_j/L)\hat f(k)\phi_k(x).
\ee
We have proved in \cite[Theorem~4.1]{mauropap}, \cite[Theorem~2.1]{frankbern} the following:
\begin{theorem}\label{summtheo}
For every $L>0$ and integer $S>\alpha$,  we have
\be\label{phikernlocest}
|\Phi_L(x,y)|\le c\frac{L^\alpha}{\max(1, (L\rho(x,y))^S)}, \qquad x,y\in\XX,
\ee
and
\be\label{phikernl1est}
\sup_{x\in\XX}\int_\XX |\Phi_L(x,y)|d\mu(y) \le c.
\ee
Consequently, for $1\le p\le \infty$,
\be\label{sigmaopbd}
\|\sigma_L(f)\|_{\mu;p}\le c\|f\|_{\mu;p}, \qquad f\in L^p.
\ee
\end{theorem}
We will also need the following two propositions. Proposition~\ref{criticalprop} is proved in \cite[Proposition~5.1]{eignet1}. The definition of regular measures in \cite{eignet1} is different from the one in this paper, but Proposition~\ref{mzequivprop} shows that they are equivalent. 

\begin{prop}\label{criticalprop}
Let $d>0$, $S>\a$ be an integer,  and \eref{ballmeasurecond}, \eref{singlegaussbd} hold. Let $\nu$ satisfy $\tn\nu\tn_{R,d}<\infty$,  $L>0$, and $\kappa_1$ be as in \eref{ballmeasurecond}. Let $1\le p\le \infty$. In this proposition, all constants will depend upon $S$.\\
{\rm (a)} If $g_1:[0,\infty)\to [0,\infty)$ is a nonincreasing function, then for any $L>0$, $r>0$, $x\in\XX$,
\be\label{g1ineq}
L^\a\int_{\Delta(x,r)}g_1(L\rho(x,y))d|\nu|(y)\le \frac{2^{\a}(\kappa_1+(d/r)^\a)\a}{1-2^{-\a}}\tn \nu\tn_{R,d}\int_{rL/2}^\infty g_1(u)u^{\a-1}du.
\ee
{\rm (b)} If  $r\ge  1/L$, then
\be\label{phiintaway}
\int_{\Delta(x,r)}|\Phi_L(x,y)|d|\nu|(y) \le c_1(1+(dL)^\a)(rL)^{-S+\a}\tn \nu\tn_{R,d}.
\ee
{\rm (c)} We have
\be\label{phiinttotal}
\int_\XX|\Phi_L(x,y)|d|\nu|(y)\le c_2(1+(dL)^\a)\tn \nu\tn_{R,d},
\ee 
\be\label{philpnorm}
\|\Phi_L(x,\circ)\|_{\nu;\XX,p} \le c_3L^{\a/p'}(1+(dL)^\a)^{1/p}\tn \nu\tn_{R,d}.
\ee
\end{prop}

In the sequel, we will assume $S>\a$ to be a fixed, large integer, and will not indicate the dependence of the constants on $S$.

Next, we recall the following Proposition~\ref{overlapprop},  proved essentially in  \cite[Eqn.~(4.40), Theorem~2.2]{frankbern}. 

\begin{prop}\label{overlapprop}
Let $L\ge 1$, $\C=\{x_1,\cdots,x_M\}\subset \XX$, $\kappa>1$, $A\ge 2$, $\delta(\C)\le \kappa q(\C)$. Let $X_k=B(x_k,\delta(\C))$, $\tilde X_k=B(x_k,A\delta(\C))$. In the following, all constants will depend upon $A$ and $\kappa$. Then for every $P\in\Pi_L$, we have
\be\label{overlapineql1wodiff}
\sum_{k=1}^M \mu(X_k)\|P\|_{\mu;\infty,\tilde X_k} \le c\{(\delta(\C) L)^\a +\min(1,(\delta(\C) L)^{\a -S})\}\|P\|_{\mu;1},
\ee
\be\label{derineqwanted}
\sum_{k=1}^M \mu(X_k)\|\tn \nabla P\tn_\circ\|_{\mu;\infty,\tilde X_k} \le cL\{(\delta(\C) L)^\a +\min(1,(\delta(\C) L)^{\a -S})\}\|P\|_{\mu;1}, 
\ee
and
\be\label{linfbernstein}
\|\tn\nabla P\tn_\circ\|_{\mu;\infty} \le cL\|P\|_{\mu;\infty}.
\ee
\end{prop}

In our proof of this proposition in  \cite{frankbern}, we used $A=2$, but the same proof works in the more general case, verbatim, except for the following changes (using equation numbers and notations from \cite{frankbern}) : The set ${\cal I}$ defined after (4.34) should be redefined by ${\cal I}= \{j : \rho(x, \tilde X_j) \ge (2A+1)\delta\}$, and the two displayed equations after (4.34) are changed to
$$
|\rho(x,y)-\rho(x,\tilde X_j)|=|\rho(x,y)-\rho(x,z_j)|\le \rho(y,z_j)\le 2A\delta\le (2A/(2A+1))\rho(x,\tilde X_j),
$$
and
$$
\delta\le (2A+1)^{-1}\rho(x,\tilde X_j)\le \rho(x,y) \le \frac{4A+1}{2A+1}\rho(x,\tilde X_j)
$$
respectively. We prefer not to reproduce the entire proof to accommodate these minor changes.
We need to prove an $L^p$ analogue of the above proposition.
\begin{lemma}\label{overlaplemmalp}
Let $L\ge 1$, $\C=\{x_1,\cdots,x_M\}\subset \XX$, $A\ge 2$, $\kappa>1$, $\delta(\C)\le \kappa q(\C)$. Let $X_k=B(x_k,\delta(\C))$, $\tilde X_k=B(x_k,A\delta(\C))$. Then for every $P\in\Pi_L$, we have
\be\label{overlapineqlpwodiff}
\left\{\sum_{k=1}^M \mu(X_k)\|P\|_{\infty,\tilde X_k}^p\right\}^{1/p} \le c\{(\delta(\C) L)^\a +\min(1,(\delta(\C) L)^{\a -S})\}^{1/p}\|P\|_{\mu;p},
\ee
and
\be\label{gradientoverlapineq}
\left\{\sum_{k=1}^M \mu(X_k)\|\tn\nabla P\tn_\circ\|_{\infty,\tilde X_k}^p\right\}^{1/p} \le cL\{(\delta(\C) L)^\a +\min(1,(\delta(\C) L)^{\a -S})\}^{1/p}\|P\|_{\mu;p}.
\ee
\end{lemma}

\begin{Proof}\ 
We will prove
\eref{gradientoverlapineq}. The proof of \eref{overlapineqlpwodiff} is similar, but simpler, and is ommitted. We observe first that for any differentiable $f :\XX\to\RR$, and $x\in \XX$,
$$
\tn\nabla f\tn_x =\sup\langle \nabla f(x), F(x)\rangle_x,
$$
where the supremum is over all  vector fields $F$ with $\tn F\tn_x=1$.

Let $F$ be an arbitrarily fixed  vector field with $\tn F\tn_x=1$ for all $x\in\XX$. We use Theorem~\ref{tensorrieszthorintheo} with $L^\infty(\mu)$ in place of $\X$, $p_0=q_0=1$, $p_1=q_1=\infty$, $w_k=\mu(X_k)$, and 
$$
{\mathcal U}f(x) =(\chi_1(x)F(x)(\sigma_{2L}(f)),\cdots,\chi_M(x)F(x)(\sigma_{2L}(f))),
$$
where each $\chi_k$ is the characteristic function of $\tilde X_k$. Then, for $1\le p<\infty$,
$$
\|{\mathcal U}f\|_{L^\infty(\mu),{\bf W},p}=\left\{\sum_{k=1}^M\mu(X_k)\|F(\sigma_{2L}(f))\|_{\mu;\infty,\tilde X_k}^p\right\}^{1/p},
$$
and the formula holds also for $p=\infty$ with an obvious modification.
Using \eref{derineqwanted}, \eref{linfbernstein}, \eref{sigmaopbd} with $2L$ in place of $L$, $\sigma_{2L}(f)$ in place of $P$, we then see that for $p=1,\infty$, $f\in L^p(\mu)$,
$$
\begin{array}{ll}
\|{\mathcal U}f\|_{L^\infty(\mu), {\bf W}, p}  &\le cL\{(\delta(\C) L)^\a +\min(1,(\delta(\C) L)^{\a -S})\}^{1/p}\|\sigma_{2L}(f)\|_{\mu;p}\\[2ex]
& \le c_1L\{(\delta(\C) L)^\a +\min(1,(\delta(\C) L)^{\a -S})\}^{1/p}\|f\|_{\mu;p}.
\end{array}
$$
Theorem~\ref{tensorrieszthorintheo} now implies
$$
\left\{\sum_{k=1}^M \mu(X_k)\|F(\sigma_{2L}(f))\|_{\mu;\infty,\tilde X_k}^p\right\}^{1/p} 
\le cL\{(\delta(\C) L)^\a +\min(1,(\delta(\C) L)^{\a -S})\}^{1/p}\|f\|_{\mu;p}.
$$
Since $F$ is an arbitrary  unit vector field, this leads to
$$
\left\{\sum_{k=1}^M \mu(X_k)\|\tn\nabla(\sigma_{2L}(f))\tn_\circ\|_{\mu;\infty,\tilde X_k}^p\right\}^{1/p} 
\le cL\{(\delta(\C) L)^\a +\min(1,(\delta(\C) L)^{\a -S})\}^{1/p}\|f\|_{\mu;p}.
$$
Since $\sigma_{2L}(P)=P$ for $P\in\Pi_L$, this implies  \eref{gradientoverlapineq}. 
\end{Proof}

\subsection{Krein's extension theorem}\label{kreinsect}
The purpose of this section is to prove the following special case of the Krein extension theorem. Let $\X$ be a normed linear space, ${\cal K}$ be a  subset of its normed dual $\X^*$, and ${\cal V}$ be a linear subspace of $\X$. We say that a linear functional $x^*\in {\cal V}^*$ is positive on ${\cal V}$ with respect to ${\cal K}$ if $x^*(f)\ge 0$ for every $f\in {\cal V}$ with the property that $y^*(f)\ge 0$ for every $y^*\in {\cal K}$.
\begin{theorem}\label{kreintheo}
Let $\X$ be a normed linear space, ${\cal K}$ be a bounded subset of its normed dual $\X^*$, ${\cal V}$ be a linear subpace of $\X$, $x^*\in {\cal V}^*$ be positive on ${\cal V}$ with respect to ${\cal K}$. We assume further that there exists $v_0\in {\cal V}$ such that $\|v_0\|_\X=1$ and
\be\label{kreinbetacond}
\inf_{y^*\in{\cal K}} y^*(v_0) =\beta^{-1}>0.
\ee
Then there exists an extension $X^*\in\X^*$ of $x^*$ which is positive on $\X$ with respect to ${\cal K}$ and satsifies 
\be\label{kreinextnormest}
\|X^*\|_{\X^*}\le \beta\sup_{y^*\in {\cal K}} \|y^*\|_{\X^*}x^*(v_0).
\ee
\end{theorem}

\begin{Proof}\ 
In this proof only, let $M=\sup_{y^*\in {\cal K}} \|y^*\|_{\X^*}$, and for $f_1,f_2\in\X$, we will say that $f_1\succeq f_2$ if $y^*(f_1)\ge y^*(f_2)$ for every $y^*\in {\cal K}$. In this proof only, let 
$$
p(f)=\inf\{ x^*(P) : P\in {\cal V},\ P\succeq f\}.
$$
For any $f\in\X$ and $y^*\in{\cal K}$, we have $|y^*(f)|\le M\|f\|_\X \le \beta M \|f\|_\X y^*(v_0)= y^*(\beta M \|f\|_\X v_0)$. Since $\pm \beta M\|f\|_\X v_0\in {\cal V}$, it follows that $p(f)$ is a finite number for $f\in\X$.
 It is not difficult to check that $p$ is a sublinear functional; i.e.,
$$
p(f_1+f_2)\le p(f_1)+p(f_2), \ p(\gamma f_1)=\gamma p(f_1), \quad f_1,f_2\in \X, \ \gamma\ge 0.
$$
If $P, Q\in {\cal V}$, and $P\succeq Q$ then the fact that $x^*$ is positive on ${\cal V}$ with respect to ${\cal K}$ implies that $x^*(P)\ge x^*(Q)$. So, $p(Q)=x^*(Q)$ for all $Q\in {\cal V}$. The Hahn--Banach theorem \cite[Theorem~(14.9), p.~212]{hewittbk} then implies that there exists an extension of $x^*$ to a linear functional $X^*$ on $\X$ such that $X^*(f)\le p(f)$, $f\in \X$. Then
$$
X^*(f)=-X^*(-f)\ge -p(-f)=\sup\{x^*(-P) : P\in {\cal V},\ P\succeq -f\} = \sup\{x^*(Q) : Q\in {\cal V}, \ f\succeq Q\}.
$$
This implies two things. First, let $f\succeq 0$. Choosing $Q$ in the last supremum expression to be $0$, we see that $X^*(f)\ge 0$. Second, as we observed earlier, $\beta M \|f\|_\X v_0 \succeq f\succeq -\beta M \|f\|_\X v_0$. Since $\pm \beta M\|f\|_\X v_0\in {\cal V}$, we obtain that $|X^*(f)| \le \beta M x^*(v_0)\|f\|_\X$. This proves \eref{kreinextnormest}, and in particular, that $X^*\in\X^*$. 
\end{Proof}

\bhag{Proofs of the results in Section~\ref{mainsect}.}\label{proofsect}

We start with the proof of Theorem~\ref{mztheolp}. 

\noindent\textsc{Proof of Theorem~\ref{mztheolp}.}  \\
We assume that $1<p<\infty$, $\C=\{x_1,\cdots,x_M\}$; the case $p=1$ is simpler, and is essentially done in \cite[Theorem~3.2]{frankbern}. We use the notation $X_k=B(x_k,\delta(\C))$, $\tilde X_k=B(x_k,A\delta(\C))$. Using the fact that $\nabla |P|^p = p|P|^{p-1}\mbox{ sgn }(P)\nabla P$, we deduce that for any $k=1,\cdots,M$, $z,y\in \tilde X_k$,
$$
\left||P(z)|^p-|P(y)|^p\right|\le 2Ap\delta(\C)\|P\|_{\infty,\tilde X_k}^{p-1}\|\tn\nabla P\tn_\circ\|_{\infty,\tilde X_k}.
$$
We may assume that $L\delta(\C)\le 1$. Hence, using H\"older's inequality, \eref{overlapineqlpwodiff}, and \eref{gradientoverlapineq}, we obtain that
\begin{eqnarray*}
\lefteqn{\sum_{k=1}^M\mu(X_k)\sup_{z,y\in \tilde X_k}\left||P(z)|^p-|P(y)|^p\right|\le 2Ap\delta(\C)\sum_{k=1}^M\mu(X_k)\|P\|_{\infty,\tilde X_k}^{p-1}\|\tn\nabla P\tn_\circ\|_{\mu;\infty,\tilde X_k}}\\
&\le&2Ap\delta(\C)\left\{\sum_{k=1}^M\mu(X_k)\|P\|_{\mu;\infty,\tilde X_k}^p\right\}^{1/p'}\left\{\sum_{k=1}^M\mu(X_k)\|\tn\nabla P\tn_\circ\|_{\mu;\infty,\tilde X_k}^p\right\}^{1/p}\\
&\le& cAp\delta(\C) L\|P\|_{\mu;p}^{p/p'}\|P\|_{\mu;p} =cAp\delta(\C) L\|P\|_{\mu;p}^p.
\end{eqnarray*}
With $c_2=1/c$, this proves \eref{mzineqstronglp} if $LAp\delta(\C)\le c_2\eta$.

To prove part (b), we observe that
\begin{eqnarray*}
\lefteqn{\left|\int_\XX |P(z)|^pd\mu(z)-\sum_{k=1}^M \mu(Y_k)|P(x_k)|^p\right|=\left|\sum_{k=1}^M \int_{Y_k} \{|P(z)|^p-|P(x_k)|^p\}d\mu(z)\right|}\\
&\le& \sum_{k=1}^M \int_{Y_k} \left||P(z)|^p-|P(x_k)|^p\right|d\mu(z)\le \sum_{k=1}^M \mu(Y_k)\sup_{z,y\in Y_k}\left||P(z)|^p-|P(y)|^p\right|\\
&\le&c \sum_{k=1}^M \mu(X_k)\sup_{z,y\in \tilde X_k}\left||P(z)|^p-|P(y)|^p\right|.
\end{eqnarray*}
Hence, \eref{discstronglpmzineq} follows from \eref{mzineqstronglp}. 

The proof of part (c) is easier.
Let $|P(z^*)|=\|P\|_{\mu;\infty}$. By definition of $\delta(\C)$, there exists $x^*\in\C$ such that $\rho(z^*,x^*)\le \delta(\C)$. Then in view of \eref{linfbernstein}, we have
$$
|P(z^*)|-|P(x^*)|\le \delta(\C)\|\tn\nabla P\tn_\circ\|_{\mu;\infty} \le cL\delta(\C)\|P\|_{\mu;\infty};
$$
i.e.,
$$
\max_{x\in\C}|P(x)|\le \|P\|_{\mu;\infty} \le |P(x^*)| + cL\delta(\C)\|P\|_{\mu;\infty}\le \max_{x\in\C}|P(x)|+cL\delta(\C)\|P\|_{\mu;\infty}.
$$
This leads to \eref{mzineqstronginf}.
\qed

In the proofs of the other results in Section~\ref{mainsect}, we will often need the following observation. If $K\subseteq\XX$ is a compact subset and $\e>0$, we will say that a subset $\C\subseteq K$ is $\e$--separated if $\rho(x,y)\ge \e$ for every $x,y\in \C$, $x\not=y$. Since $K$ is compact, there exists a finite, maximal $\e$--separated subset $\{x_1,\cdots,x_M\}$ of $K$. If $x\in K\setminus \cup_{k=1}^M B(x_k,\e)$, then $\{x,x_1,\cdots,x_M\}$ is a strictly larger $\e$--separated subset of $K$. So, $K\subseteq 
\cup_{k=1}^M B(x_k,\e)$. Moreover, the balls $B(x_k,\e/2)$ are mutually disjoint. 

The proof of Theorem~\ref{genpartitiontheo} requires some further preparation. First, we recall a lemma \cite[Lemma~4.4]{frankbern}. For a set $Y$, we denote the cardinality of $Y$ by $|Y|$.
\begin{lemma}\label{noofinterlemma}
Let $\C$ be a finite set for which \eref{cuniformitygen} holds, $\delta(\C)\le 2\kappa$, $A>0$, and $x\in \XX$.  Then 
$$
\left|\{y\in\C : x\in B(y,A\delta(\C))\}\right|\le c_1(1+A)^{2\a},
$$
 where $c_1>0$ is independent of $A$ and $\delta(\C)$. In particular,
$$
\left|\{y\in\C : B(x,A\delta(\C)\cap B(y, A\delta(\C))\not=\emptyset\}\right| \le c_1(1+2A)^{2\a}.
$$
\end{lemma}

\begin{Proof} 
In this proof, let $\delta:=\delta(\C)$.  Let $y_1,\cdots,y_m\in \C$ and $x\in \cap_{k=1}^m B(y_k,A\delta)$. Then $B(x,\delta)\subseteq\cap_{k=1}^m B(y_k,(1+A)\delta)$. Since $q(\C)\ge \delta/\kappa$, the balls $B(y_k,\delta/(2\kappa))$ are pairwise disjoint, and their union is a subset of $B(x,(1+A)\delta)$. Therefore, \eref{doublingcond} implies that
 \begin{eqnarray*}
\mu(B(x,\delta))&\le& \min_{1\le k\le m}\mu(B(y_k,(1+A)\delta))\le \frac{1}{m}\sum_{k=1}^m\mu(B(y_k,(1+A)\delta))\le \frac{c(1+A)^\a}{m}\sum_{k=1}^m\mu(B(y_k,\delta/(2\kappa)))\\
&=&\frac{c(1+A)^\a}{m}\mu\left(\cup_{k=1}^m B(y_k,\delta/(2\kappa))\right) \le \frac{c(1+A)^\a}{m}\mu(B(x,(1+A)\delta))\le \frac{c_1(1+A)^{2\a}}{m}\mu(B(x,\delta)).
\end{eqnarray*}
Thus, $m\le c_1(1+A)^{2\a}$. 
\end{Proof}

 The following lemma is needed in the construction of the partition in Theorem~\ref{newstrongmztheo}. The proof is based on some ideas in the book  \cite[Appendix 1]{davidbk} of David.
\begin{lemma}\label{partitionlemma}
Let $\tau$ be a positive measure on $\XX$, ${\cal A}$ be a finite subset of $\XX$ satisfying
$$
q({\cal A})/2 \le \delta({\cal A}) \le \kappa q({\cal A})
$$
for some $\kappa>0$, $\{Z_y\}_{y\in {\cal A}}$ be a partition of $\XX$ such that each $Z_y\subseteq B(y,\gamma\delta({\cal A}))$ for some $\gamma\ge 1$. (We do not require that each $Z_y$ be nonempty.) Then there exists a subset ${\cal G}\subseteq {\cal A}$ and a partition $\{Y_y\}_{y\in {\cal G}}$ such that $Z_y\subseteq Y_y$, $\tau(Y_y)\ge  c\min_{z\in{\cal A}} \tau(B(z,\gamma\delta({\cal A})))$, and each $Y_y\subseteq B(y,3\gamma\delta({\cal A}))$. In particular, $\delta({\cal G})\le 3\gamma\delta({\cal A})$ and $q({\cal G})\ge q({\cal A})$.
\end{lemma}
\begin{Proof}\ 
In this proof, let $\delta=\delta({\cal A})$, $m=\min_{z\in{\cal A}} \tau(B(z,\gamma\delta))$. In view of Lemma~\ref{noofinterlemma}, at most $c^{-1}$ of the balls $B(y,\gamma\delta)$ can intersect each other, the constant $c$ depending upon $\gamma$ and $\kappa$.  Let ${\cal G}=\{y\in{\cal A} : \tau(Z_y)\ge cm\}$. Now, we define a function $\phi$ as follows. If $z\in {\cal G}$, we write $\phi(z)=z$. Otherwise, let $z\in {\cal A}\setminus {\cal G}$. Since $\{Z_y\}$ is a partition of $\XX$, we have
$$
m\le \tau(B(z, \gamma\delta)) =\sum_{y\in {\cal A}}\tau(B(z, \gamma\delta)\cap Z_y).
$$
Since each $Z_y\subseteq B(y,\gamma\delta)$, it follows that at most $c^{-1}$ of the $Z_y$'s have a nonempty intersection with  $B(z,\gamma\delta)$. So, there must exist $y\in {\cal A}$ for which
$$
\tau(B(z, \gamma\delta)\cap Z_y) \ge cm.
$$
Clearly, each such $y\in {\cal G}$. We imagine an enumeration of ${\cal A}$, and among the $y$'s for which $\tau(B(z, \gamma\delta)\cap Z_y)$ is maximum, pick the one with the lowest index. We then define $\phi(z)$ to be this $y$. Necessarily,
$\phi(z)=y\in {\cal G}$, and $B(z, \gamma\delta)\cap Z_y\subseteq B(z, \gamma\delta)\cap B(y,\gamma\delta)$ is nonempty. So, 
\be\label{pf3eqn1}
\rho(z,\phi(z))\le 2\gamma\delta, \quad B(z, \gamma\delta)\subseteq B(\phi(z),3\gamma\delta), \quad \tau(B(z, \gamma\delta)\cap Z_{\phi(z)}) \ge cm.
\ee
Now, we define
$$
Y_y= \cup\{Z_z : \phi(z)=y, z\in {\cal A}\}, \qquad y\in {\cal G}.
$$
 For each $z\in {\cal A}$, $Z_z\subseteq Y_{\phi(z)}$. Since $Z_z$ is a partition of $\XX$, $\XX=\cup_{y\in {\cal G}} Y_y$. If $x\in \XX$, $x\in Y_y\cap Y_{y'}$ for $y, y'\in {\cal G}$, then $x\in Z_z$ with $\phi(z)=y$ and $x\in Z_{z'}$ with $\phi(z')=y'$. Since $\{Z_z\}$ is a partition of $\XX$, it follows that $z=z'$, and hence $y=y'$. Thus, $\{Y_y\}$ is a partition of $\XX$, $\tau(Y_y)\ge \tau(Z_y)\ge cm$, and 
$$
Y_y\subseteq \cup_{\phi(z)=y}Z_z \subseteq \cup_{\phi(z)=y}B(z, \gamma\delta)\subseteq B(y,3\gamma\delta).
$$
\end{Proof}

\noindent
\textsc{Proof of Theorem~\ref{genpartitiontheo}.} The partition $Y_k$ is required to satisfy three goals: firstly, we wish to ensure that $\mu(Y_k)\sim d^\a$, secondly, we wish to be able to obtain a lower bound on $|\nu|(Y_k)$ as stated in Theorem~\ref{genpartitiontheo}, and finally, we wish to ensure that $x_k\in Y_k$. Accordingly, we will start with an appropriately dense subset of $\supp(\nu)$, and construct a corresponding partition in a somewhat obvious manner. We will then use Lemma~\ref{partitionlemma} three times to ensure the three goals; first with $\mu$ in place of $\tau$, then with $|\nu|$ in place of $\tau$ with the resulting partition, and finally build a somewhat artificial measure $\tau$ supported on the finite subset obtained in the second step, and use the lemma with this measure. This  ensures that each of the sets in the resulting partition contains at least one point of the set in the second step, but not necessarily the set created in the third step. However, this is easy to arrange with an increase in the mesh norm by a constant factor.

In this proof only, let ${\cal G}_1=\{y_1,\cdots,y_N\}$ be a maximal $d/2$--separated subset of $\supp(\nu)$. Then $\supp(\nu)\subseteq \cup_{k=1}^N B(y_k,d/2)$, and $d/4\le q({\cal G}_1)/2\le \delta({\cal G}_1)\le 3d/2$. We let $Z_{y_1,1}=B(y_1,\delta({\cal G}_1))$, and for $k=2,\cdots,N$, $Z_{y_k,1}=B(y_k,\delta({\cal G}_1))\setminus \cup_{j=1}^{k-1} B(y_j,\delta({\cal G}_1))$. Then $\{Z_{y,1}\}_{y\in {\cal G}_1}$ is a partition of $\XX$ and each $Z_{y,1}\subseteq B(y,\delta({\cal G}_1))$, $y\in {\cal G}_1$. 

We apply Lemma~\ref{partitionlemma} first with $\mu$ in place of $\tau$, resulting in a subset ${\cal G}_2\subseteq{\cal G}_1$, $\delta({\cal G}_2)\le 3\delta({\cal G}_1)$, and a partition $\{Z_{y,2}\}_{y\in {\cal G}_2}$ of $\XX$ such that for each  $y\in {\cal G}_2$, $Z_{y,2}\subseteq B(y,3\delta({\cal G}_1))\subseteq B(y,3\delta({\cal G}_2))$,  and $c_1\min_{z\in{\cal G}_1}\mu(B(z,\delta({\cal G}_1)))\le \mu(Z_{y,2}) \le \mu(B(y,3\delta({\cal G}_1)))$; i.e., $\mu(Z_{y,2})\sim (\delta({\cal G}_1))^\a\sim d^\a$. 

We apply Lemma~\ref{partitionlemma} again with ${\cal G}_2$ in place of ${\cal A}$, $\{Z_{y,2}\}$ as the corresponding partition, and $|\nu|$ in place of $\tau$. This yields a subset ${\cal G}_3\subseteq{\cal G}_2$ and a partition $\{Z_{y,3}\}_{y\in {\cal G}_3}$ of $\XX$ with $\delta({\cal G}_3)\le 3\delta({\cal G}_2)$,  such that for each each $y\in {\cal G}_3$, $Z_{y,2}\subseteq Z_{y,3}\subseteq B(y,3\delta({\cal G}_2))\subseteq B(y,9\delta({\cal G}_1))\cap B(y,3\delta({\cal G}_3))$, and $|\nu|(Z_{y,3})\ge c_2\min_{z\in {\cal G}_2}|\nu|(B(z,\delta({\cal G}_2)))=:u$. Since ${\cal G}_2\subseteq \supp(\nu)$, $u$ is a positive number. We note that $\mu(Z_{y,3})\sim d^\a$ as well. 

At this point, we still have not proved that  $Z_{y,3}\cap {\cal G}_3$ is nonempty for each $y\in {\cal G}_3$. Towards this end, we repeat an application of Lemma~\ref{partitionlemma} with the measure, to be denoted in this proof only by $\tau$, that associates the mass $u>0$ with each $y\in {\cal G}_3$. This gives us a subset ${\cal G}_4 \subseteq {\cal G}_3$ and a partition $\{Z_{y,4}\}_{y\in{\cal G}_4}$ with $\delta({\cal G}_4)\le 3\delta({\cal G}_3)$, such that for each $y\in {\cal G}_4$, $Z_{y,2}\subseteq Z_{y,3}\subseteq Z_{y,4}\subseteq B(y, 3\delta({\cal G}_3))\subseteq B(y, 27\delta({\cal G}_1))$, $\mu(Z_{y,4})\sim d^\a$ and $\tau(Z_{y,4})\ge c_3\min_{z\in {\cal G}_3}\tau(B(z, \delta({\cal G}_3)))\ge c_3u>0$. Necessarily, each $Z_{y,4}$ contains some element of ${\cal G}_3$. 

We pick one element from each $Z_{y,4}\cap {\cal G}_3$ to form the set $\C=\{x_1,\cdots,x_M\}$, and rename the set $Z_{y,4}$ containing $x_k$ to be $Y_k$. By construction, $\{Y_k\}_{k=1}^M$ is a partition of $\XX$ with each $x_k\in Y_k \subseteq B(x_k, 54\delta({\cal G}_1))$, $\mu(Y_k)\sim d^\a$ and $|\nu|(Y_k)\ge c_3u$. 
\qed

We are now in a position to prove Theorem~\ref{newstrongmztheo}.

\noindent\textsc{Proof of Theorem~\ref{newstrongmztheo}.} 
We  apply Theorem~\ref{mztheolp} with the set $\C$ and the partition $\{Y_k\}$. We observe that for any $k$,
\begin{eqnarray*}
\lefteqn{\left|\int_{Y_k}|P(y)|^pd\mu(y)-\frac{\mu(Y_k)}{|\nu|(Y_k)}\int_{Y_k}|P(z)|^pd|\nu|(z)\right|}\\
&=&\left|\int_{Y_k}\left\{|P(y)|^p- \frac{1}{|\nu|(Y_k)}\int_{Y_k}|P(z)|^pd|\nu|(z)\right\}d\mu(y)\right|\\
&=&\frac{1}{|\nu|(Y_k)}\left|\int_{Y_k}\int_{Y_k}\left\{|P(y)|^p-|P(z)|^p\right\}d|\nu|(z)d\mu(y)\right|\\
&\le& \mu(Y_k)\max_{z,y\in Y_k}\left||P(y)|^p-|P(z)|^p\right|.
\end{eqnarray*}
Since $Y_k\subseteq B(x_k,81d)$, \eref{mzineqstronglp} leads to \eref{newstrongmzineq}. This completes the proof of part (a). Part (b) follows from Theorem~\ref{mztheolp}(c).
\qed

We find it convenient to prove Proposition~\ref{mzequivprop} next, so that we may use such statements as Proposition~\ref{criticalprop} which were proved with the definition as in \eref{regreconcile} rather than the one which have used in this paper.

\noindent
\textsc{Proof of Proposition~\ref{mzequivprop}.} In the proof of part (a) only, let $\lambda>\tn\nu\tn_{R,d}$, $r>0$, $x\in\XX$, and let $\{y_1,\cdots,y_N\}$ be a maximal $2d/3$--separated subset of $B(x,r+2d/3)$. Then $B(x,r)\subseteq B(x,r+2d/3)\subseteq \cup_{j=1}^N B(y_j,2d/3)$. So,
$$
|\nu|(B(x,r)) \le |\nu|(B(x,r+2d/3))\le \sum_{j=1}^N |\nu|(B(y_j,2d/3))\le \sum_{j=1}^N |\nu|(B(y_j,d)) \le \lambda Nd^\a.
$$
The balls $B(y_j,d/3)$ are mutually disjoint, and $\cup_{j=1}^N B(y_j,d/3)\subseteq B(x,r+d)$. In view of \eref{ballmeasurelowbd}, $d^\a\le c\mu(B(y_j,d/3))$ for each $j$. So,
$$
|\nu|(B(x,r))\le \lambda Nd^\a\le c\lambda\sum_{j=1}^N \mu(B(y_j,d/3))=c\lambda\mu(\cup_{j=1}^N B(y_j,d/3))\le c\lambda \mu(B(x,r+d)).
$$
Since $\lambda>\tn\nu\tn_{R,d}$ was arbitrary, this leads to the first inequality in \eref{regreconcile}. The second inequality follows from \eref{ballmeasurecond}. The converse statement is obvious. This completes the proof of part (a).

The second estimate in \eref{regequiv} is clear from the definitions. The first estimate in \eref{regequiv} follows by applying \eref{regreconcile} with $r=\gamma d$.
\qed

In the proof of Theorem~\ref{mzmeascharacttheo}, we will often need the following observation.
\begin{lemma}\label{philowbdlemma}
There exists a constant $\beta\in (0,1/2)$ such that for any $L>0$, $x\in\XX$, 
\be\label{phillowbd}
|\Phi_L(x,y)|\ge (1/2)\Phi_L(x,x)\ge cL^\a, \qquad \rho(x,y)\le \beta/L, \quad \
\ee
Hence, for $1\le p<\infty$ and $\nu\in {\cal M}$,
\be\label{phillplocbd}
\int_{B(x,\beta/L)}|\Phi_L(x,y)|^pd|\nu| \ge cL^{\a p}|\nu|(B(x,\beta/L)).
\ee
\end{lemma}
\begin{Proof}\ 
In this proof only, let $P(y)=\Phi_L(x,y)$, $y\in \XX$, and $|P(y^*)|=\|P\|_{\mu;\infty}$. Then $P(x)=\Phi_L(x,x)$, and Schwarz inequality and \eref{christbd} show that 
$$
\|P\|_{\mu;\infty}=|P(y^*)|\le \Phi_L(x,x)^{1/2}\Phi_L(y^*,y^*)^{1/2}\le c_1L^\a \le c_2\Phi_L(x,x)=c_2P(x)\le c_2\|P\|_{\mu;\infty}.
$$
Thus,  $\|P\|_{\mu;\infty} \sim L^\a$. Since $P\in\Pi_L$, we conclude from \eref{linfbernstein} that 
$$
|P(y)-P(x)|\le c_3L\rho(x,y)\|P\|_{\mu;\infty}=c_4L\rho(x,y)P(x).
$$
Hence, with $\beta=\min(1/2,1/(2c_4))$, we obtain that $|P(y)| \ge (1/2)P(x) \ge cL^\a$ if $\rho(x,y)\le \beta/L$. 
\end{Proof}

We are now in a position to prove Theorem~\ref{mzmeascharacttheo}

\noindent
\textsc{Proof of Theorem~\ref{mzmeascharacttheo}.} Let $\nu$ be $1/L$--regular, and without loss of generality, $\|\nu\|_{R,1/L}=1$. Using \eref{phiinttotal} with $\mu$ and $\nu$,  and the fact that both $\mu$ and $\nu$ are $1/L$--regular, we deduce that
$$
\sup_{x\in\XX}\max\left\{\int_\XX |\Phi_{2L}(x,y)|d\mu(y), \int_\XX |\Phi_{2L}(x,y)|d|\nu|(y)\right\}\le c.
$$
Therefore, using Fubini's theorem, we conclude as in \cite[Corollary~5.2]{eignet1} that for $p=1,\infty$,
$$
\|\sigma_{2L}(f)\|_{\nu;p}\le c\|f\|_{\mu;p}.
$$
Hence, the Riesz--Thorin interpolation theorem shows that this inequality is valid also for all $p$, $1\le p\le \infty$. If $P\in \Pi_L$, we use this inequality with $P$ in place of $f$, and recall that $\sigma_{2L}(P)=P$ to deduce that $\|P\|_{\nu;p}\le c\|P\|_{\mu;p}$, as claimed in the first part of Theorem~\ref{mzmeascharacttheo}(a).

Conversely, suppose that for some $p$, $1\le p<\infty$,
\be\label{pf4eqn1}
\|P\|_{\nu;p}^p\le A\|P\|_{\mu;p}^p, \qquad P\in \Pi_L.
\ee 
Let $x\in\XX$. We apply \eref{pf4eqn1} with $P(y)=\Phi_L(x,y)$. Using \eref{philpnorm} with the $1/L$--regular measure $\mu$ in place of $\nu$, we see that $\|P\|_{\mu;p}^p \le cL^{\a p/p'}=cL^{\a(p-1)}$.  Therefore, \eref{pf4eqn1} and \eref{phillplocbd} together imply that with $\beta$ as defined in Lemma~\ref{philowbdlemma},
$$
c_1L^{\a p}|\nu|(B(x,\beta/L)) \le \|P\|_{\nu;p}^p\le cAL^{\a(p-1)}.
$$
Thus, $\nu$ is $\beta/L$--regular, with $\tn\nu\tn_{R,\beta/L}\le (c/c_1)A$. In view of \eref{regequiv} applied with $1/\beta$ in place of $\gamma$, this completes the proof of the converse statement in Theorem~\ref{mzmeascharacttheo}(a).

The proof of the first part of Theorem~\ref{mzmeascharacttheo}(b) relies upon Theorem~\ref{newstrongmztheo}(a). Let $1\le p<\infty$. With the constants $c_1$, $c_2$ as in that theorem, let $c_4=c_2/(8p)$, $L\ge 4c_4\max(81, c_1^{-1})$, $d=4c_4/L$. Then $d\le c_1$, and $\nu$ is $d/4$--dominant. From the definition, this means that for every $x\in\XX$, $\nu(B(x,d/4))>0$. In particular, each $B(x,d/4)\cap \supp(\nu)$ is nonempty; i.e., $\delta(\supp(\nu))\le d/4<d \le 1/81$. Therefore, all the conditions of Theorem~\ref{newstrongmztheo}(a) are satisfied so that 
  \eref{newstrongmzineq} holds with $\eta=1/2$. This shows that for every such $p$,
\be\label{pf5eqn1}
(1/2) \|P\|_{\mu;p}^p \le \sum_{k=1}^M\frac{\mu(Y_k)}{|\nu|(Y_k)}\int_{Y_k}|P(z)|^pd|\nu|(z).
\ee
Here, we recall that $\mu(Y_k)\sim d^\a$, and $|\nu|(Y_k)\ge c\min_{x\in\XX}|\nu|(B(x,d/4))$.  Since $\nu$ is $d/4$--dominant, we have $|\nu|(Y_k)\ge c\tn\nu\tn_{D,d/4}^{-1}d^\a$. So, \eref{pf5eqn1} leads to
$$
\|P\|_{\mu;p}^p \le c_3\tn\nu\tn_{D,d/4}\sum_{k=1}^M\int_{Y_k}|P(z)|^pd|\nu|(z) = c_3\tn\nu\tn_{D,c_4/L}\|P\|_{\nu;p}^p.
$$
This completes the proof of the first part of Theorem~\ref{mzmeascharacttheo}(b).

Finally, we prove the converse statement in Theorem~\ref{mzmeascharacttheo}(b). Accordingly, we assume that $\nu$ is $1/L$--regular. In view of part (a) of this theorem and the assumption of the converse statement, we have with $A\sim \tn\nu\tn_{R,1/L}^{-1/p}$,
\be\label{pf5eqn2}
A\|P\|_{\nu;p}^p \le \|P\|_{\mu;p}^p \le A_1\|P\|_{\nu;p}^p, \qquad P\in\Pi_L.
\ee
We will use Lemma~\ref{philowbdlemma} as before. Let $x\in\XX$, and $P(y)=\Phi_L(x,y)$, and $r\ge 1$ to be chosen later. 
Using \eref{phillplocbd} with $\mu$ in place of $\nu$, and \eref{ballmeasurelowbd}, we obtain
\be\label{pf5eqn3}
\frac{c_7}{A_1}L^{\a(p-1)}\le A_1^{-1}\|P\|_{\mu;p}^p \le \|P\|_{\nu;p}^p =\int_{B(x,r/L)} |P(y)|^pd|\nu|(y) + \int_{\Delta(x,r/L)}|P(y)|^pd|\nu|(y).
\ee
Since $\nu$ is assumed to be $1/L$--regular, we may apply \eref{phiintaway} to conclude that
\be\label{pf5eqn4}
\int_{\Delta(x,r/L)}|P(y)|^pd|\nu|(y)\le \|P\|_{\mu;\infty}^{p-1}\int_{\Delta(x,r/L)}|P(y)|d|\nu|(y) \le c_8L^{\a(p-1)}r^{\a-S}\tn\nu\tn_{R,1/L}.
\ee
We now choose 
$$
r=\max\left(1, \left(\frac{2c_8A_1\tn\nu\tn_{R,1/L}}{c_7}\right)^{1/(S-\a)}\right).
$$
 Then \eref{pf5eqn3} and \eref{pf5eqn4} together lead to
$$
\frac{c_7}{2A_1}L^{\a(p-1)}\le \int_{B(x,r/L)} |P(y)|^pd|\nu|(y) \le c_9L^{\a p}|\nu|(B(x,r/L)).
$$
Thus,
$$
|\nu|(B(x,r/L))\ge c_{10}A_1^{-1}L^{-\a}.
$$
This completes the proof of the converse statement.
\qed

\noindent
\textsc{Proof of Theorem~\ref{quadtheo}.}
We find a finite subset $\C=\{x_1,\cdots,x_M\}$ and a partition $\{Y_k\}$ as in Theorem~\ref{genpartitiontheo}.  In view of \eref{subsetdelta},  each $x_k\in Y_k\subseteq B(x_k,324\delta(\C))$. Moreover, the conditions of Theorem~\ref{mztheolp} are satisfied with appropriate $\kappa$ for this $\C$. In view of  \eref{doublingcond} and \eref{mzineqstronglp}, we conclude that for a suitably chosen $c$, $L\le c\delta(\C)^{-1}\sim d^{-1}$, 
\be\label{pf6eqn5}
\sum_{k=1}^M \mu(B(x_k, 324\delta(\C)))\max_{y,z\in B(x,324\delta(\C))}|P(y)-P(z)|\le (1/4)\|P\|_{\mu;1}.
\ee
 In this proof only, let 
$$
x_k^*(P)=\frac{1}{|\nu|(Y_k)}\int_{Y_k}P(z)d|\nu|(z), \qquad k=1,\cdots,M.
$$ 
Then for $P\in\Pi_L$, and $k=1,\cdots,M$,
\bea\label{pf6eqn1}
\lefteqn{\left|\int_{Y_k}|P(y)|d\mu(y)-\mu(Y_k)|x_k^*(P)|\right|}\nonumber\\
&=&\left|\int_{Y_k}\left\{|P(y)|- |x_k^*(P)|\right\}d\mu(y)\right|\nonumber\\
&\le& \int_{Y_k}\left||P(y)|- |x_k^*(P)|\right|d\mu(y)\nonumber\\
&\le& \int_{Y_k}|P(y)- x_k^*(P)|d\mu(y)\nonumber\\
&\le& \frac{1}{|\nu|(Y_k)}\int_{Y_k}\int_{Y_k}|P(y)-P(z)|d|\nu|(z)d\mu(y)\nonumber\\
&\le& \mu(B(x,324\delta(\C)))\max_{y,z\in B(x,324\delta(\C))}|P(y)-P(z)|.
\eea

Then \eref{pf6eqn1} and \eref{pf6eqn5} imply that
$$
\left|\sum_{k=1}^M \int_{Y_k}|P(y)|d\mu(y) -\sum_{k=1}^M \mu(Y_k)|x_k^*(P)|\right|\le (1/4)\|P\|_{\mu;1};
$$
i.e.,
\be\label{pf6eqn2}
(3/4)\|P\|_{\mu;1}\le \sum_{k=1}^M \mu(Y_k)|x_k^*(P)|\le (5/4)\|P\|_{\mu;1}.
\ee

Moreover, if each $x_k^*(P)\ge 0$, then the same estimate as \eref{pf6eqn1} with $P$ in place of $|P|$ leads to
$$
\left|\int_\XX P(y)d\mu(y)-\sum_{k=1}^M \mu(Y_k)x_k^*(P)\right|\le (1/4)\|P\|_{\mu;1}\le (1/3)\sum_{k=1}^M \mu(Y_k)x_k^*(P).
$$
Thus, if each $x_k^*(P)\ge 0$, then 
\be\label{pf6eqn4}
\int_\XX P(y)d\mu(y)\ge (2/3)\sum_{k=1}^M \mu(Y_k)x_k^*(P)\ge 0. 
\ee

Now, we wish to use Theorem~\ref{kreintheo}. We let $\X$ be the space $\RR^M$, equipped with the norm $\|{\bf y}\|=\sum_{k=1}^M\mu(Y_k)|y_k|$, where ${\bf y}=(y_1,\cdots,y_k)$. For the set ${\cal K}$, we choose the set of coordinate functionals; $y_k^*({\bf y})=y_k$. Then ${\cal K}$ is clearly a compact subset of $\X^*$. We consider the operator ${\cal S} :\Pi_L\to\RR^M$ given by $P\mapsto (x_1^*(P),\cdots,x_M^*(P))$, and take the subspace ${\cal V}$ of $\X$ to be the range of ${\cal S}$. The lower estimate in \eref{pf6eqn2} shows that ${\cal S}$ is invertible on ${\cal V}$. We define the functional $x^*$ on ${\cal V}$ by 
$$
x^*({\cal S}(P))=\int_\XX P(z)d\mu(z)-(1/3)\sum_{k=1}^M \mu(Y_k)x_k^*(P), \qquad P\in \Pi_L.
$$
 Our observations in the previous paragraph show that $x^*$ is positive on ${\cal V}$ with respect to ${\cal K}$. The element $(1,\cdots,1)\in {\cal V}$ serves as $v_0$ in Theorem~\ref{kreintheo}. Theorem~\ref{kreintheo} then implies that there exists a nonnegative functional $X^*$ on ${\cal X}=\RR^M$ that extends $x^*$. We may identify this functional with $(\tilde W_1,\cdots,\tilde W_M)\in\RR^M$, such that each $\tilde W_k\ge 0$. The fact that $X^*$ extends $x^*$ means that for each $P\in\Pi_L$,
\be\label{pf6eqn3}
\int_\XX P(x)d\mu(x)=\sum_{k=1}^M (\tilde W_k+(1/3)\mu(Y_k))x_k^*(P)=:\sum_{k=1}^M \frac{W_k}{|\nu|(Y_k)}\int_{Y_k}P(y)d|\nu|(y).
\ee
Writing $W(y)=W_k/|\nu|(Y_k)$ for $y\in Y_k$, we have now proved \eref{diffpolyquad}. By construction, $W_k\ge (1/3)\mu(Y_k)\ge c_1d^\a$.  If $\nu$ is $d$--regular, then $|\nu|(Y_k)\le |\nu|(B(x_k,81d))\le c_2\tn\nu\tn_{R,d}d^\a$. Hence, $W(y)\ge c\tn\nu\tn_{R,d}^{-1}$ for all $y\in\XX$.
\qed

The proof of Theorem~\ref{posmztheo} uses the following lemma proved in  \cite[Lemma~5.5]{eignet1}:
\begin{lemma}\label{quaderrlemma}
Let the product assumption hold, and $L>0$. If $\nu$ is a quadrature measure of order $2A^*L$, $|\nu|(\XX)\le c$, and $P_1,P_2\in\Pi_{2L}$ then for any $p, r$, $1\le p, r\le\infty$ and any positive number $R>0$,
\be\label{polyquaderr}
\left|\int_\XX P_1P_2d\mu-\int_\XX P_1P_2d\nu\right| \le  c_1L^{2\a}\e_L\|P_1\|_{\mu;p}\|P_2\|_{\mu;r} \le c(R)L^{-R}\|P_1\|_{\mu;p}\|P_2\|_{\mu;r}.
\ee
\end{lemma}

\noindent
\textsc{Proof of Theorem~\ref{posmztheo}.}  
 Let $x\in\XX$, and $P:=\Phi_L(x,\circ)\in \Pi_L$. Taking $\beta$ as in Lemma~\ref{philowbdlemma}, we obtain from \eref{phillplocbd} applied with $\tau $ in place of $\nu$ that 
$$
c_1L^{2\a}\tau(B(x,\beta/L))\le \int_{B(x,\beta/L)} |P(y)|^2d\tau(y)\le \int_\XX  |P(y)|^2d\tau(y).
$$
Since $\tau$ is a positive quadrature measure of order $2A^*L$, we now obtain from Lemma~\ref{quaderrlemma} used with $P_1=P_2=P$, $p=r =1$, $R= 1$ that
\begin{eqnarray*}
c_1L^{2\a}\tau(B(x,\beta/L)) \le  \int_\XX |P(y)|^2d\mu(y) + c\|P\|_{\mu;1}^2/L= \Phi_L(x,x) +c_2\|P\|_{\mu;1}^2/L.
\end{eqnarray*}
In view of \eref{phikernl1est} and  \eref{christbd}, $\|P\|_{\mu;1}\le c_3$,  $\Phi_L(x,x)\le c_3L^\a$. We deduce that for sufficiently large $L$:
$$
L^{2\a}\tau(B(x,\beta/L))\le c_4L^\a;
$$
i.e., $\tau(B(x,\beta/L))\le c_4L^{-\a}$. In view of Theorem~\ref{mzmeascharacttheo}, this implies that $\|P\|_{\tau;p}\le c_5\|P\|_{\mu;p}$ for all $p$ with $1\le p <\infty$.

It remains to prove that\be\label{pf7eqn2}
\|P\|_{\mu;p}\le c\|P\|_{\tau;p}, \qquad P\in \Pi_L, \quad 1\le p<\infty.
\ee
Towards this end, we introduce the discretized version of the operator $\sigma_L$:
$$
\sigma_L(\tau;f,x)=\int_\XX \Phi_L(x,y)f(y)d\tau(y), \qquad f\in L^1(\tau), \ x\in\XX, \ L>0.
$$
If $f\in L^p(\tau)$ and $g\in L^{p'}(\mu)$, then it is easy to verify using Fubini's theorem that
\begin{eqnarray*}
\lefteqn{\int_\XX \sigma_L(\tau;f,x)g(x)d\mu(x)=\int_\XX\int_\XX \Phi_L(x,y)f(y)d\tau(y)g(x)d\mu(x)}\\
&=& \int_\XX\int_\XX \Phi_L(x,y)g(x)d\mu(x)f(y)d\tau(y)=\int_\XX\sigma_L(g,y)f(y)d\tau(y).
\end{eqnarray*}
Using the duality principle and \eref{sigmaopbd}, we conclude that
\begin{eqnarray*}
\|\sigma_L(\tau;f)\|_{\mu;p}&=&\sup_{\|g\|_{\mu;p'}=1}\left|\int_\XX \sigma_L(\tau;f,x)g(x)d\mu(x)\right| =\sup_{\|g\|_{\mu;p'}=1}\left|\int_\XX\sigma_L(g,y)f(y)d\tau(y)\right|\\
&\le&\sup_{\|g\|_{\mu;p'}=1}\|\sigma_L(g)\|_{\mu;p'}\|f\|_{\tau;p} \le c\|f\|_{\tau;p}.
\end{eqnarray*}
In particular,
\be\label{pf7eqn3}
\|\sigma_{2L}(\tau;f)\|_{\mu;p}\le c\|f\|_{\tau;p}.
\ee
Since we do not assume that the product of polynomials in $\Pi_{2L}$ is itself in $\Pi_{2A^*L}$, it does not follow from the fact that $\tau$ is a quadrature measure of order $2A^*L$ that $\sigma_{2L}(\tau;P)=P$ for $P\in \Pi_L$. Nevertheless, Lemma~\ref{quaderrlemma} allows us to conclude that $\sigma_{2L}(\tau;P)$ is close to $P$ if $P\in\Pi_L$. Let $x\in \XX$. Since $\tau$ is a quadrature measure of order $2A^*L$,  we may use Lemma~\ref{quaderrlemma} with $R=1$, $r=1$, $P$ in place of $P_1$ and $\Phi_{2L}(x,\circ)$ in place of $P_2$ to conclude that
\begin{eqnarray*}
\left| \sigma_{2L}(\tau;P,x)-P(x)\right|&=&\left| \sigma_{2L}(\tau;P,x)- \sigma_{2L}(P,x)\right|\\
&=&\left|\int_\XX P(y)\Phi_{2L}(x,y)d\tau(y)-\int_\XX P(y)\Phi_{2L}(x,y)d\mu(y)\right|\\
&\le& cL^{-1}\|\Phi_{2L}(x,\circ)\|_{\mu;1}\|P\|_{\mu;p}.
\end{eqnarray*}

In view of \eref{phikernl1est}, and the fact that $\mu$ is a probability measure, this implies that
$$
\left|\|\sigma_{2L}(\tau;P)\|_{\mu;p}-\|P\|_{\mu;p}\right|\le \|\sigma_{2L}(\tau;P)-P\|_{\mu;p}\le cL^{-1}\|P\|_{\mu;p}.
$$
Thus, if $L\ge 2c$, then \eref{pf7eqn3} implies
$$
\|P\|_{\mu;p}\le 2\|\sigma_{2L}(\tau;P)\|_{\mu;p}\le c_4\|P\|_{\tau;p}.
$$
This proves \eref{pf7eqn2} and completes the proof of the theorem. \qed

\renewcommand{\theequation}{\Alph{section}.\arabic{equation}}
\setcounter{section}{0}
\setcounter{equation}{0}
\renewcommand{\thesection}{\Alph{section}}
\section{Appendix}

In this appendix, we follow an idea in the paper \cite{pesenson} of Geller and Pesenson to prove that the product assumption is valid when $\phi_k$ (respectively, $\ell_j^2$) are eigenfunctions (respectively, eigenvalues) of selfadjoint, uniformly elliptic partial differential operators of second order satisfying some technical conditions. In the appendix, we will use the symbol $q$ again to denote the dimension of the manifold $\XX$. 

First, we need to recall the notion of exponential maps on the manifold. For any $x\in \XX$, there exists a neighborhood $V$ of $x$, a number $\epsilon>0$ and a $C^\infty$ mapping $\gamma :(-2,2)\times {\cal U} \to \XX$, where ${\cal U}=\{(y,w)\in T\XX; y\in V, w\in T_y\XX, \tn w\tn_y<\e\}$, 
such that $\gamma(\circ,y,w)$ is the unique geodesic of $\XX$ with $\gamma(0,y,w)=y$, and the tangent vector at $y$ being $w$ (\cite[Proposition~2.7]{docormo2}). 
In this appendix, we will denote by $B_\e(0)$ the open Euclidean ball in $\RR^q$ with center at $0$ and radius $\e$. For any tangent space $T_x\XX$, we may consider an appropriate coordinate chart, and view $B_\e(0)$ as a subset of $T_x\XX$, with $0$ corresponding to $x$. The \emph{exponential map} at $x$ is the mapping $\exp_x : B_\e(0) \subset T_x\XX\to \XX$ defined by $\exp_x(w)=\gamma(1,x,w)$, where $\gamma$ is the mapping just described. Thus, $\exp_x(w)$ is the point on $\XX$ where one reaches by following the geodesic at $x$, with tangent vector given by $w/\|w\|$ for a length of $\|w\|$. For every $x\in \XX$, there exists an $\e>0$ such that $\exp_x :B_\e(0)\subset T_x\XX\to \XX$ is a diffeomorphism of $B_\e(0)$ onto an open subset of $\XX$.  Since $\XX$ is compact, we may choose $\e$ to be the same for all $x\in\XX$. The largest value of such $\e$ is called the \emph{injectivity radius} of $\XX$, to be denoted in this appendix by $\iota$.  If $\delta>0$ and $\overline{B_\delta(0)}\subset B_\iota(0)\subset T_x\XX$, then $\exp_x(B_\delta(0))$ is called a \emph{normal ball} of radius $\delta$ centered at $x$. Normal neighborhoods of $x$ are defined in the obvious way. If $\{\partial_j\}$ is a basis for $T_x\XX$, the \emph{normal coordinate system} (with respect to $\{\partial_j\}$) at $x$ is defined on a normal neighborhood of $x$ by $\disp\x(u_1,\cdots,u_q)=\exp_x\left(\sum_{j=1}^q u_j\partial_j\right)$.

Let $\mathcal{P}$ be a self-adjoint differential operator of second order. In terms of a normal coordinate system at a point $x\in\XX$ the operator $\mathcal{P}$ can be expressed in the form 
$$
{\mathcal P}f=\sum_{\k\in\ZZ^q,\ |\k|\le 2} a_{\k,x}(\u)\frac{\partial^{|\k|} f}{\partial \u^\k}.
$$
 The operator is \emph{strongly uniformly elliptic} if there are constants $c_1, c_2>0$ (independent of $x$) such that
$$
c_1\|\y\|^2 \le \sum_{|\k|=2}a_{\k,x}(\u)\y^\k \le c_2\|\y\|^2, \qquad \u\in B_\iota(0), \ x\in\XX,\ \y\in\RR^q.
$$
We assume that there exists a constant $C>1$ such that
\be\label{pdecond1}
\left|\frac{\partial^{|\m|} a_{\k,x}}{\partial \u^\m}\right|\le C^{|\m|}, \qquad x\in\XX, \ \u\in B_{\iota}(0),\ \m\in\ZZ^q.
\ee
The eigenvalues of ${\mathcal P}$ can be enumerated in the form $\{\ell_k^2\}$ ($\ell_k\uparrow \infty$), and we let $\phi_k$ be the eigenfuction corresponding to $\ell_k^2$. We assume (by choosing a larger $C$ if necessary) that
\be\label{pdecond2}
\sum_{\ell_k\le L}\left|\frac{\partial^{|\m|}\phi_k(x)}{\partial \u^\m}\right|^2\le C^{2|\m|}L^{q+2|\m|}, \qquad x\in\XX,\ \m\in\ZZ^q, \ L>0.
\ee
This result follows essentially from the estimates on the derivatives of the heat kernel corresponding to ${\mathcal P}$ given by
Kordyukov \cite[Theorem~5.5]{kordyukov91}, and the Tauberian theorem in our paper \cite[Proposition~4.1]{frankbern}, except for the dependence of the constants involved on $\m$.  An immediate consequence of \eref{pdecond2} is the following. If $Q=\sum_k \hat Q(k)\phi_k\in\Pi_L$, then
$$
\left|\frac{\partial^{|\m|}Q(x)}{\partial {\bf u}^\m}\right|^2=\left|\sum_k\hat Q(k)\frac{\partial^{|\m|}\phi_k(x)}{\partial {\bf u}^\m}\right|^2\le \left\{\sum_k |\hat Q(k)|^2\right\}\left\{\sum_{\ell_k\le L}\left|\frac{\partial^{|\m|}\phi_k(x)}{\partial {\bf u}^\m}\right|^2\right\} \le C^{2|\m|}L^{q+2|\m|}\|Q\|_{\mu;2}^2;
$$
i.e.,
\be\label{mockbernstein}
\left|\frac{\partial^{|\m|}Q(x)}{\partial {\bf u}^\m}\right|\le L^{q/2}(CL)^{|\m|}\|Q\|_{\mu;2}.
\ee

We are now ready to prove the product assumption, in fact, a much stronger statement:
\begin{theorem}\label{prodtheo}
Let ${\mathcal P}$ be a second order, strongly uniformly elliptic, self-adjoint, partial differential operator on a smooth, compact manifold $\XX$ (without boundary), the eigenvalues of ${\mathcal P}$  be enumerated in the form $\{\ell_k^2\}$ ($\ell_k\uparrow \infty$), and we let $\phi_k$ be the eigenfuction corresponding to $\ell_k^2$. Assume further that \eref{pdecond1} and \eref{pdecond2} are satisfied. There exists $A^*\ge 2$ such that if $Q, R\in \Pi_L$, then $QR\in \Pi_{A^*L}$.
\end{theorem}

\begin{Proof}\ 
Let $N\ge 2$ be an integer. In view of Leibniz's formula, one can write 
$$
{\mathcal P}^N(QR)(x)=\sum_{\k,\m\in\ZZ^q,\ |\k|+|\m|\le 2N} b_{\k,\m}(x)  \frac{\partial^{|\m|}Q(x)}{\partial {\bf u}^\m}\frac{\partial^{|\k|}R(x)}{\partial {\bf u}^\k},
$$
where $b_{\k,\m}$'s are products of derivatives of the coefficients $a_{\k,x}$ in ${\mathcal P}$. In view of \eref{pdecond1} and \eref{pdecond2}, we conclude that for some $A^*\ge 2$,
\be\label{pf8eqn1}
\|{\mathcal P}^N(QR)\|_{\mu;2}\le \|{\mathcal P}^N(QR)\|_{\mu;\infty}\le L^q(A^*L/2)^{2N}\|Q\|_{\mu;2}\|R\|_{\mu;2}.
\ee
In this proof only, let for $f\in L^2(\mu)$,
$$
S_L (f)=\sum_{\ell_k\le A^*L}\hat f(k)\phi_k.
$$
We observe that $\|f-S_L (f)\|_{\mu;2}\le \|f\|_{\mu;2}$.
Since ${\mathcal P}\phi_k=\ell_k^2\phi_k$, it follows that ${\mathcal P}^N(S_L(QR))\!=\!S_L({\mathcal P}^N(QR))$. Consequently,
\bea\label{pf8eqn2}
\left\|{\mathcal P}^N(QR-S_L(QR))\right\|_{\mu;2}&=&\left\|{\mathcal P}^N(QR)-S_L({\mathcal P}^N(QR))\right\|_{\mu;2}\le \|{\mathcal P}^N(QR)\|_{\mu;2}\nonumber\\
&\le& L^q(A^*L/2)^{2N}\|Q\|_{\mu;2}\|R\|_{\mu;2}.
\eea
On the other hand,  Parseval's identity shows that
\begin{eqnarray*}
\left\|{\mathcal P}^N(QR-S_L(QR))\right\|_{\mu;2}^2&=&\sum_{\ell_k>A^*L}\ell_k^{4N}\widehat{(QR)}(k)^2\ge (A^*L)^{4N}\sum_{\ell_k>A^*L}\widehat{(QR)}(k)^2\\
&\ge&  (A^*L)^{4N}\|QR-S_L(QR)\|_{\mu;2}^2.
\end{eqnarray*}
Together with \eref{pf8eqn2}, this implies that for every  integer $N\ge 2$,
$$
\|QR-S_L(QR)\|_{\mu;2}\le L^q\left(\frac{A^*L}{2A^*L}\right)^{2N}\|Q\|_{\mu;2}\|R\|_{\mu;2}=L^q4^{-N}\|Q\|_{\mu;2}\|R\|_{\mu;2}.
$$
Letting $N\to \infty$, we conclude that $QR=S_L(QR)\in \Pi_{A^*L}$. This completes the proof.
\end{Proof}



\begin{thebibliography}{99}
\def\vol#1{{\bf #1}} 
\def\sandarbh#1 #2,, #3,, #4..{\bibitem{#1} {\sc #2}\ {\it #3}\ #4}
\def\pustak#1{{\rm ``#1'',\/}}

\sandarbh{boothby} W. M. Boothby,,, \pustak{An Introduction to Differentiable Manifolds and Riemannian Geometry},, Academic Press/Elsevier, Amsterdam, 2003...


\sandarbh{docormo1} M. P. do Carmo,,, \pustak{Differential Geometry of Curves and Surfaces},, Prentice Hall, New Jersey, 1976...

\sandarbh{docormo2} M. P. do Carmo,,, \pustak{Riemannian Geometry},, Birkh\"auser, Boston, 1992...

\sandarbh{chengliyau81} S. Y. Cheng, P. Li and S.-T. Yau,,, On the upper estimate of the heat kernel of a complete Riemannian manifold,,, Amer. J. Math., {\bf 103} (1981) no.5, 1021--1063...

\sandarbh{achaspissue} C. K. Chui and D. L. Donoho,,, Special Issue: Diffusion maps and wavelets,,, Applied and Computational Harmonic Analysis, vol. 21(1), 2006...
\sandarbh{mauro1} R. R. Coifman and M. Maggioni,,, Diffusion wavelets,,, Appl. Comput. Harmon. Anal. \vol{21} (2006), 53--94... 
\sandarbh{damelin} S. B. Damelin and J. Levesley,,, Hyperinterpolation, packing and spherical
designs on projective spaces,,, Manuscript, private communication...
\sandarbh{davidbk} G. David,,, \pustak{Wavelets and singular integrals on curves and surfaces},, Lecture notes in mathematics, \vol{1465}, Springer Verlag, Berlin, 1992...

\sandarbh{davies97} E. B. Davies,,, $L^p$ spectral theory of higher-order elliptic differential operators,,, Bull.
London Math. Soc., {\bf 29} (1997) 513--546...
\sandarbh{dungey06} N. Dungey,,, Some remarks on gradient estimates for heat kernels,,,  Abstr. Appl. Anal.  2006, Art. ID 73020, 10 pp...

\sandarbh{frankbern} F. Filbir and H. N. Mhaskar,,, A quadrature formula for diffusion polynomials corresponding to a generalized heat kernel,,, Accepted for publication in J. Four. Anal. Appl...

\sandarbh{pesenson} D. Geller and I. Z. Pesenson,,, Band--limited localized Parseval frames and Besov spaces on compact homogeneous manifolds,,, arXiv:1002.3841v1 [math.FA], 2010...
\sandarbh{grigoryan95} A. Grigor'yan,,, Upper bounds of derivatives of the heat kernel on an arbitrary complete manifold,,, Journal of Functional Analysis {\bf 127} (1995), no. 2, 363--389...
\sandarbh{grigoryan99} A. Grigor'yan,,, Estimates of heat kernels on Riemannian manifolds,,, Spectral Theory and Geometry (Edinburgh,
1998) (E. B. Davies and Y. Safarov, eds.), London Math. Soc. Lecture Note Ser., vol. 273,
Cambridge University Press, Cambridge, 1999, pp. 140--225...
\sandarbh{wtheatkernel} A. Grigor'yan,,,  Heat kernels on weighted manifolds and applications,,,  Cont. Math. 398 (2006) 93--191...
\sandarbh{grigoryanheatmetric} A. Grigor'yan,,, Heat kernels and function theory on metric
measure spaces,,, to appear in "Handbook of Geometric Analysis No.2" ed. L. Ji, P. Li, R. Schoen, L. Simon, Advanced Lectures in Math., IP,   2008...
\sandarbh{besov} K. Hesse, H. N. Mhaskar, and I. H. Sloan,,, Quadrature in Besov spaces on the Euclidean sphere,,,  Journal of Complexity, {\bf 23} (2007), 528--552...
\sandarbh{sloansharma} K. Hesse and I. H. Sloan,,, Hyperinterpolation on the sphere,,, in  ``Frontiers in interpolation and approximation'',   Pure Appl. Math. (Boca Raton), 282, Chapman \& Hall/CRC, Boca Raton, FL, 2007, pp. 213--248...
\sandarbh{hewittbk} E. Hewitt and K. Stromberg,,, \pustak{Real and abstract analysis},, Springer, New York, 1975...
\sandarbh{horn} R. A. Horn and C. R. Johnson,,, \pustak{Matrix analysis},, Cambridge University Press, 1985...
\sandarbh{jms} P. W. Jones, M. Maggioni, R. Schul,,, Universal local parametrizations via heat kernels and eigenfunctions of the Laplacian,,, Manuscript arXiv:0709.1975...
\sandarbh{kordyukov91} Yu. A. Kordyukov,,, Lp-theory of elliptic differential operators on manifolds of bounded geometry,,, Acta Appl. Math. {\bf 23} (1991), no. 3, 223--260...

\sandarbh{lafon} S. Lafon,,, Diffusion maps and geometric harmonics,,, PhD thesis, Yale University, Dept of Mathematics \& Applied Mathematics, 2004...
\sandarbh{lee} J. M. Lee,,, \pustak{Introduction to smooth manifolds},, Springer, New York, 2003...

\sandarbh{quadconst} Q. T. Le Gia and H. N. Mhaskar,,, Localized linear polynomial operators and  quadrature formulas on the sphere,,, SIAM J. Numer. Anal. {\bf 47} (1) (2008), 440--466...

\sandarbh{lubinskymz} D. S. Lubinsky,,, Marcinkiewicz inequalities: methods and results,,, in Recent Progress in
Inequalities (G. V. Milovanovic, Ed.), pp. 213--240, Kluwer, Dordrecht, 1998...
\sandarbh{mauropap} M. Maggioni and H. N. Mhaskar,,, Diffusion polynomial frames on 
metric measure spaces,,, Appl.  Comput. Harm. Anal., \vol{24} (3) (2008), 329--353...

\sandarbh{eignet1} H. N. Mhaskar,,, Eignets for function approximation on manifolds,,, To appear in Appl. Comput. Harm. Anal...
\sandarbh{mnw1} H. N. Mhaskar, F. J. Narcowich and J. D. Ward,,, Spherical Marcinkiewicz-Zygmund inequalities and positive quadrature,,, Math. Comp. {\bf 70} (2001), no. 235, 1113--1130.   (Corrigendum: Math. Comp. {\bf 71} (2001), 453--454.)..
\sandarbh{oneil} B. Oneill,,, \pustak{Elementary differential geometry},, Elsevier, Amsterdam, 2006...


\sandarbh{taylorbk} M. E. Taylor,,, \pustak{Partial differential equations, basic theory},, Springer Verlag, New York, 1996...

\sandarbh{trigub} R. M. Trigub and E. S. Belinsky,,, \pustak{Fourier analysis and approximation of functions},, Kluwer, Dodrecht, 2004...
\sandarbh{zygmund} A. Zygmund,,, {\rm ``Trigonometric
Series'',\/},, Cambridge University Press, Cambridge, 1977...

\end{thebibliography}
\end{document}